\nonstopmode
\input amstex
\input xy
\xyoption{all}
\newdir{ >}{{}*!/-10pt/@{>}}
\message{ ) }
\input amsppt.sty
\hsize 30pc
\vsize 47pc
\def\nmb#1#2{#2}         
\def\cit#1#2{\ifx#1!\cite{#2}\else#2\fi} 
\def\totoc{}             
\def\ign#1{}             
\redefine\o{\circ}

\define\al{\alpha}

\define\ka{\kappa}

\define\si{\sigma}
\define\ta{\tau}
\define\ph{\varphi}

\define\De{\Delta}

\define\CC{\Bbb{C}}
\define\oo{\infty}
\redefine\i{^{-1}}
\define\x{\times}
\let\on=\operatorname
\def\sr#1%
{\ifmmode{}^\dagger\else${}^\dagger$\fi\ifvmode
\vbox to 0pt{\vss
 \hbox to 0pt{\hskip\hsize\hskip1em
 \vbox{\hsize3cm\eightpoint\raggedright\pretolerance10000
 \noindent #1\hfill}\hss}\vss}\else
 \vadjust{\vbox to0pt{\vss%
 \hbox to 0pt{\hskip\hsize\hskip1em%
 \vbox{\hsize3cm\eightpoint\raggedright\pretolerance10000%
 \noindent #1\hfill}\hss}\vss}}\fi%
}
\def\today{\ifcase\month\or
 January\or February\or March\or April\or May\or June\or
 July\or August\or September\or October\or November\or December\fi
 \space\number\day, \number\year}
\def\p{\partial}
\NoRunningHeads
\topmatter
\title
Lifting mappings over invariants of finite groups
\endtitle
\author Andreas Kriegl, Mark Losik, Peter W. Michor, Armin Rainer
\endauthor
\leftheadtext{M.\ Losik, P. W.\ Michor}
\address
A\. Kriegl: Institut f\"ur Mathematik, Universit\"at Wien,
Strudlhofgasse 4, A-1090 Wien, Austria
\endaddress
\email Andreas.Kriegl\@univie.ac.at \endemail
\address
M. Losik: Saratov State University, ul. Astrakhanskaya, 83,
410026 Saratov, Russia
\endaddress
\email LosikMV\@info.sgu.ru \endemail
\address
P\. W\. Michor: Institut f\"ur Mathematik, Universit\"at Wien,
Strudlhofgasse 4, A-1090 Wien, Austria; {\it and:}
Erwin Schr\"odinger Institut f\"ur Mathematische Physik,
Boltzmanngasse 9, A-1090 Wien, Austria
\endaddress
\email Peter.Michor\@esi.ac.at \endemail
\address
A\. Rainer: Institut f\"ur Mathematik, Universit\"at Wien,
Strudlhofgasse 4, A-1090 Wien, Austria
\endaddress
\email armin\_rainer\@gmx.net\endemail
\abstract
We characterize those regular, holomorphic or formal maps into the orbit
space $V/G$ of a complex representation of a finite group $G$ which admit a
regular, holomorphic or formal lift to the representation space $V$. In
particular, the case of complex reflection groups is investigated.
\endabstract

\thanks{M.L., P.W.M., and A.R.\ were supported
     by `Fonds zur F\"orderung der
     wissenschaftlichen For\-schung,
     Projekt P~14195~MAT'.}
\endthanks
\endtopmatter
\document

\head\totoc\nmb0{1}. Introduction \endhead
The problem of lifting of morphisms to the orbit space for real or
complex representations of compact Lie groups was studied in several papers.

The existence of lifts of smooth diffeomorphisms of orbit spaces for a real
representation of a compact Lie group $G$ in a real vector space $V$
was investigated in \cit!{2}, \cit!{12}, and \cit!{6}. The condition
for the existence of lifts of
smooth and analytic curves in orbit spaces for real representations of
compact Lie groups was found in \cit!{1}. Note that in the real case
the orbit space $V/G$ has a natural structure of a real semialgebraic subset
of some space $\Bbb R^m$ and its stratification coincides with the isotropy
type stratification of $V/G$. This stratification plays the main role
in the results mentioned above.

In the complex case the lifting problem is more complicated and only
the case of a representation of a finite group $G$ in a complex vector
space $V$ was studied. In this case the orbit space $Z=V/G$ coincides with
the categorical quotient $V//G$ which is a normal affine variety.
Therefore the orbit space $Z$ has the natural structure of a complex analytic
set and there are several types of morphisms into $V//G$, like
regular, rational or holomorphic. To formulate the condition of lifting
one needs to use the isotropy type (Luna) stratification of the orbit space
$V/G$ which is finer than its stratification as affine variety.
The conditions of lifting for holomorphic automorphisms of orbit spaces
were found for the Weyl groups in \cit!{9} and for any finite groups
in \cit!{4}. In \cit!{7} it was proved that each holomorphic lift of
a regular automorphism of the orbit space is regular.

In this paper we consider the conditions for lifts of germs of holomorphic
morphisms at $0$ from $\Bbb C^p$ to $Z$, for lifts of regular maps from
$\Bbb C^p$ to $Z$, and for lifts of formal morphisms from $\Bbb C^p$ to
$Z$, i.e\., the morphisms of the $\Bbb C$-algebra $\Bbb C[Z]$ to the ring of 
formal power series
$\Bbb C[[X_1,\dots,X_p]]$ in variables $X_1,\dots,X_p$.
These conditions are formulated with the use of the spaces $J^q_0(\Bbb C^p,Z)$
and $J^\infty_0(\Bbb C^p,Z)$ of jets at $0\in\Bbb C^p$ of orders $q$ and
$\infty$, respectively. In particular,
we consider these conditions in the case when $G$ is a finite group
generated by complex reflections. Finally, we calculate the above conditions
for some simple cases.

Algebraically, the above problems of lifting could be reformulated as partial
cases of the general problem of extensions for morphisms of the category
of $\Bbb C$-algebras from a subalgebra to the whole algebra.

Note that, by Lefschetz's principle (see, for example, \cit!{14})
the conditions of lifting which are formulated in algebraic terms
are valid for any field of characteristic zero.

In section \nmb!{2} we formulate the types of lifting problems which are
solved in the paper and introduce the tools we need for this: Luna's
stratification and the jet spaces for affine varieties, in particular,
for a $G$-module $V$ and for the orbit space $Z=V/G$.

In section \nmb!{3} we define the functions $\tilde T(A_1,\dots,A_d)$ and
$\tilde\frak T(A_1,\dots,A_d)$, which are used in the conditions of lifting
for holomorphic, regular, and formal lifts and study their properties.

In section \nmb!{4} we obtain the conditions for local and global holomorphic lifts and
regular lifts.

In section \nmb!{5} we obtain the condition for formal lifts.

In section \nmb!{6} we consider the above conditions of lifting for complex
reflection groups and calculate some of these conditions for a reflection
group in $\Bbb C$ and for the dihedral groups.

\head\totoc\nmb0{2}. Preliminaries\endhead

\subhead{\nmb.{2.1}. Luna stratification of orbit spaces}\endsubhead
Let $V$ be an $n$-dimensional complex vector space, $G$ a finite subgroup
of $GL(V)$, and $\Bbb C[V]^G$ the algebra of $G$-invariant polynomials on
$V$.

The following facts are well known (see, for example \cit!{11}).
Denote by $Z$ the categorical quotient $V//G$, i.e\. the normal affine
algebraic variety with the coordinate ring $\Bbb C[V]^G$. Since the group
$G$ is finite, the categorical quotient $V//G$ is the geometric one,
i.e\. $V//G$ is the orbit space $V/G$.
Let $\pi=\pi_V:V\to Z$ be the quotient projection. The affine
algebraic variety $Z$ has the natural structure of a complex analytic
space: Let $\si_1,\dots,\si_m$ be the minimal system of homogeneous
generators of the algebra $\Bbb C[V]^G$ and let
$\si:V\to\si(V)\subseteq\Bbb C^m$ be the corresponding morphism.
Then $\si(V)$ is an irreducible Zariski-closed subset of $\Bbb C^m$ which
is isomorphic to the affine variety $Z$. For this presentation of $Z$
the morphism $\si:V\to\si(V)$ coincides with the projection $\pi$.

In the sequel we assume that the minimal system of homogeneous
generators $\si_1,\dots,\si_m$ is fixed and they
are called the basic generators.

Let $K$ be a subgroup of $G$. We denote by
$V_{(K)}$ the set of points of $V$ whose isotropy groups are conjugate to $K$.
By definition, $V_{(K)}\subseteq\cup_{g\in G}V^{gKg\i}$, where $V^K$ is
the subspace of $V$ of fixed points of the action of $K$ on $V$. 
Put $Z_{(K)}:=\pi(V_{(K)})$.
It is known that $\{Z_{(K)}: K < G\}$ is a finite stratification of
$Z$ into locally closed irreducible smooth algebraic subvarieties.
This is the simplest case of a Luna stratification, see \cit!{8}.
Put $V_0:=V_{(K)}$ for $K=\{\operatorname{id}\}$ and $Z_0:=\pi(V_0)$.

Denote by $Z_{>i}$ the union of the strata of codimension greater than $i$
and put $Z_{\le i}:=Z\setminus Z_{>i}$. Then $Z_{>i}$ is a Zariski-closed subset of
$Z$ and $Z_0=Z_{\le 0}=Z\setminus Z_{>0}$ is a stratum of $Z$ called
the principal stratum. Points in $Z_0$ and in $V_0$ are called regular
points.
The following proposition is evident.
\proclaim{Proposition} $Z_0$ is a Zariski-open smooth
subvariety of $Z$ and the restriction of $\pi$ to $V_0$ is an \'etale
morphism onto $Z_0$.
\endproclaim

\subhead{\nmb.{2.2}. Polarizations}\endsubhead
Let $\al$ be a homogeneous $G$-invariant polynomial of degree $d$ on $V$
and let $\al^s:V^d\to\Bbb C$ be the corresponding symmetric $d$-linear form on $V$. For
$v_1,\dots,v_k\in V$ we have
$$\align
\al(t_1v_1+\dots+t_kv_k)=&\alpha^s(t_1v_1+\dots+t_kv_k,\dots,t_1v_1+\dots+t_kv_k)\\
=&\sum_{i_1,\dots,i_d} t_{i_1}\dots t_{i_d}\;\al^s(v_{i_1},\dots,v_{i_d})\\
=&\sum_{r_1+\dots+r_k=d}t_1^{r_1}\dots t_k^{r_k}\;
\al_{r_1,\dots,r_k}(v_1,\dots,v_k), \\
\intertext{where} \\
\al_{r_1,\dots,r_k}(v_1,\dots,v_k):=&\frac{d!}{r_1!\dots r_k!}
\;\al^s(\undersetbrace \text{$r_1$ times} \to {v_1,\dots,v_1},\dots,
\undersetbrace \text{$r_k$ times}\to {v_k,\dots,v_k}) 
\endalign$$
and $t_1,\dots,t_k$ are variables.

The polynomials $\al_{r_1,\dots,r_k}$ are the usual polarizations of
$\al$ on $V^k$.

\subhead{\nmb.{2.3}. Invariant coordinates on $V$}\endsubhead
For each regular point $z_0\in Z_0$ there is
a system of regular functions $z_1,\dots,z_n$ on $Z$ such that
each $y_i:=z_i\o\pi$ equals one of the generators $\si_j$ and the functions
$z_i-z_i(z_0)$
are local parameters at $z_0$.
Then the $y_i$ are local coordinates on $V$ in a neighborhood of each point
$v\in\pi^{-1}(z_0)$. By definition, the functions $y_i$ are $G$-invariant.
These coordinates $y_i$ are called {\it invariant coordinates} on $V$.

Since we fixed the basic generators $\si_1,\dots,\si_m$,
there are only finitely many choices of such invariant coordinates on $V$.

Let $e_i$ $(i=1,\dots,n)$ be a basis of $V$ and $u_i$ the corresponding
coordinates on $V$.
Denote by $J\in\Bbb C[V]$ the jacobian $\on{det}(\frac{\partial y_i}{\partial u_j})$.
It is clear that $J$ is a homogeneous polynomial.

\proclaim{Proposition} For each integer $k>0$ there is a
$G$-invariant polynomial $\De_k\in\Bbb C[V]^G$ of minimal degree such that
$J^k$ divides $\De_k$ and the sets of zeros of $J$ and $\De_k$ coincide.
The polynomial $\De_k$ is unique up to a nonzero factor $c\in\Bbb C$.
\endproclaim

\demo{Proof} Let $J=f_1^{n_1}\dots f_s^{n_s}$ be a decomposition of $J$ into
the product of linearly independent irreducible polynomials
$f_l\in\Bbb C[V]$. Consider the principal effective divisor 
$(J)=n_1(f_1)+\dots+n_s(f_s)$ of the polynomial $J$ on $V$. 
Since for each $g\in G$ we have $J\o g=\det(g^i_j).J$, where $(g^i_j)$ is 
the matrix of $g$ in the basis $e_i$, the divisor $(J)$ is $G$-invariant. 
Then each $g\in G$ permutes the simple divisors $(f_p)$ $(p=1,\dots,s)$. 
This implies that, if $g(f_p)=(f_q)$, the coefficients $n_p$ and $n_q$ 
of the divisor $(J)$ are equal. Let $\{m_1,\dots,m_l\}$ be the set of 
distinct coefficients of the divisor $(J)$ and let, for each $m_{\al}$, 
$\varPhi_{\al}$ be the product of distinct factors $f_p$ of $J$ having 
the same power $m_{\al}$ in the above decomposition of $J$. Then we 
have $J=\prod_{\al=1}^l\varPhi_{\al}^{m_{\al}}$. By the above arguments, 
for each $\al=1,\dots,l$ the divisor $(\varPhi_{\alpha})$ of the polynomial 
$\varPhi_{\al}$ is $G$-invariant. 

Since the group $G$ is finite, for each $\al=1,\dots,l$ there is a minimal 
integer $p_{\al}>0$ such that the polynomial $\varPhi_{\al}^{p_{\al}}$ 
is $G$-invariant. 
For $\al=1,\dots,l$ let $km_{\al}=s_{\al}p_{\al}+r_{\al}$, where $s_{\al}$ 
and $r_{\al}$ are unique
nonnegative integers such that $0\le r_{\al}<p_{\al}$. Then we have
$J^k=\prod_{\al=1}^l\varPhi_{\al}^{s_{\al}p_{\al}+r_{\al}}$.
Let $\mu_{\al}$ be the least common multiple of $r_{\al}$ and $p_{\al}$. Then
$\De_k=\prod_{\al=1}^l\varPhi_{\al}^{s_{\al}p_{\al}+\mu_{\al}}$ is a
$G$-invariant polynomial of minimal degree such that the sets of zeros
of $J$ and $\De_k$ coincide and $J^k$ divides $\De_k$.

By the above formula for $J^k$, for each $G$-invariant polynomial $P$
such that the sets of zeros of $J$ and $P$ coincide and $J^k$ divides $P$,
$\De_k$ divides $P$.
\qed\enddemo

We denote by $\tilde\De_k$ the regular function on $Z$ such that
$\tilde\De_k\o\pi=\Delta_k$.
By definition, we have $\tilde\Delta_k(z_0)\neq 0$ for each $k$. Conversely,
let $y_i$ be invariant coordinates on $V$, let $z_i$ be the
corresponding regular functions on $Z$, and for some positive integer $k$
let $\tilde\Delta_k$ be the corresponding regular function on $Z$. If, for a point
$z\in Z$, we have $\tilde\Delta_k(z)\ne 0$, then $z\in Z_0$.

Later, for the sake of simplicity, we put $\De:=\De_1$ and
$\tilde\De:=\tilde\De_1$.

Denote by $V(\tilde\De)$ the set of zeros of $\tilde\De$. Thus $Z_{>0}$
is the intersection of the Zariski-closed subsets $V(\tilde\De)$ obtained
from all choices of invariant coordinates constructed from the
basic generators of $\Bbb C[V]^G$. The similar statement
is true for $V\setminus V_0$ if we replace $\tilde\De$ by $\De$.

\subhead{\nmb.{2.4}. Jet spaces}\endsubhead
Now, for an affine variety $X$, we will define the space $J^q_0(\Bbb C^p,X)$
of $q$-jets at $0$ of morphisms from $\Bbb C^p$ to $X$.
It is similar to the corresponding notion of jet spaces
for affine schemes of finite type over $\Bbb C$ in the case $p=1$
(see, for example, \cit!{10} and \cit!{5}).

We will consider now the category of $\Bbb C$-algebras. 

Let $\Bbb C[X_1,\dots,X_p]$ be the $\Bbb C$-algebra of polynomials in variables 
$X_1,\dots,X_p$ with complex coefficients and let $\frak m_p$ be the
ideal of $\Bbb C[X_1,\dots,X_p]$ generated by the $X_1,\dots,X_p$. Put
$\frak m_p^q:=(\frak m_p)^q$. Then 
$$
\frak J^q_p:=\Bbb C[X_1,\dots,X_p]/\frak m_p^{q+1}
$$
is the truncated ring of polynomials, the model jet algebra. In particular,
$\frak J_p^0=\Bbb C[X_1,\dots,X_p]/\frak m_p=\Bbb C$.

Let $A=(a_1,\dots,a_s)$ for $a_1,\dots,a_s\in\{1,\dots,p\}$ be a (unordered)
multi-index of order $|A|:=s$.
In particular, for $s=0$ we put $A:=\emptyset$.
Denote by $\frak A_{p,q}$ the set of multi-indices $A$ of
orders $\le q$. By definition, each $P\in\frak J_p^q$
can be written as $P=\sum_{A\in\frak A_{p,q}}p_A\, X_A$, where
$p_A\in\Bbb C$ and $X_A:=X_{a_1}\dots X_{a_s}$. The natural bijection 
$P\mapsto(p_A)_{A\in\frak A_{p,q}}$ between $\frak J_p^q$ 
and $\Bbb C^{\frak A_{p,q}}$ is an isomorphism of vector spaces and 
defines a structure of affine space on 
$\frak J_p^q$.  
For $q\leq r$, consider the natural morphism
$\rho_{r,q}:\frak J_p^r \to\frak J_p^q$.

For an affine variety $X$ over $\Bbb C$, the set of 
$\frak J_p^q$-valued points of $X$,
i.e\. morphisms from the coordinate ring $\Bbb C[X]$ of $X$ to the ring  
$\frak J_p^q$, is called the {\it space of $q$-jets} 
of morphisms from $\Bbb C^p$ to $X$ at $0\in\Bbb C^p$ and is denoted by 
$J^q_0(\Bbb C^p,X)$. In particular, we have
$J^0_0(\Bbb C^p,X)=X$ and $J^1_0(\Bbb C^p,X)=TX$, the total tangent bundle
of $X$.

It is evident that each polynomial function on $\frak J_p^q$ 
defines a function on $J^q_0(\Bbb C^p,X)$ and these 
functions generate a ring of $\Bbb C$-valued functions on $J^q_0(\Bbb C^p,X)$. 
It is clear that this ring is a finitely generated $\Bbb C$-algebra. 
Then $J^q_0(\Bbb C^p,X)$ supplied with this ring has a structure of an affine 
variety (not necessarily irreducible). 
For two affine varieties $X_1$ and $X_2$ 
and for a morphism $\varphi:X_1\to X_2$, we have the natural morphism  
$J^q_0(\Bbb C^p,\varphi):\bar J^q_0(\Bbb C^p,X_1)\to \bar J^q_0(\Bbb C^p,X_2)$ of 
affine varieties. Thus one can consider $J^q_0(\Bbb C^p,\quad)$ as a 
covariant functor from the category of affine varieties to itself.   

For each $h\in J^q_0(\Bbb C^p,X)$ there is a unique point $x\in X$ such that
the corresponding maximal ideal $\frak m_x$ coincides with
the kernel of the composition $\rho_{q,0}\o h$. Then the morphism $h$ can be
extended uniquely to a morphism from 
$\Cal O_x$ to $\frak J^q_p$ vanishing on
$\frak m_x^{q+1}$ and hence induces a morphism
$h_{x,q}:\Cal O_x/\frak m_x^{q+1}\to\frak J^q_p$
which, in turn, determines the initial morphism $h$ uniquely.
Therefore one can
view $J^q_0(\Bbb C^p,X)$ as the set of morphisms from the local 
rings $\Cal O_x/\frak m_x^{q+1}$ to $\frak J^q_p$
for all $x\in X$.

$$
\xymatrix{
  &\CC[X]\ar@{->}[0,1]^{h} \ar@{->}[1,0]  &\frak J^q_p \\
\frak m_x^{q+1}\ar@{^(->}[0,1]  &\Cal O_x\ar@{-->}[-1,1]^h 
\ar@{->}[0,1]  &\Cal O_x/\frak m_x^{q+1}\ar@{-->}[-1,0]_
{h_{x,q}}
}
$$
Assume that $X$ is presented as a Zariski-closed subset of $\Bbb C^m$ 
defined by
an ideal $(\varPhi_1,\dots,\varPhi_r)$ of the ring $\Bbb C[W_1,\dots,W_m]$
of polynomials with complex coefficients in variables $W_1,\dots,W_m$.

A morphism $h:\Bbb C[X]\to \frak J_p^q$ is defined by a morphism
$h':\Bbb C[W_1,\dots,W_m]\to \frak J_p^q$ with 
$h'(\varPhi_l)=0$ for $l=1,\dots,r$. 
It is determined by 
$h'(W_i)=\sum_{A\in \frak A_{p,q}}W_{i,A}\,X_A$, where $i=1,\dots,m$ and
$W_{i,A}\in\Bbb C$.
The condition $h'(\varPhi_l)=0$ is equivalent to the
vanishing of all the coefficients of the variables $X_A$ in $h(\varPhi_l)$.
Thus the map $J^q_0(\Bbb C^p,X)\to(\Bbb C^m)^{\frak A_{p,q}}$ given by
$h\mapsto (W_{i,A})_{i=1,\dots,m,\,A\in \frak A_{p,q}}$ induces a bijective
correspondence between $J^q_0(\Bbb C^p,X)$ and the Zariski-closed subset of
$(\Bbb C^m)^{\frak A_{p,q}}$ defined by $r|\frak A_{p,q}|$ many polynomial
equations, where $|\frak A_{p,q}|$ denotes the cardinality of the set $\frak A_{p,q}$.
By definition, this correspondence is an isomorphism of affine varieties.

The homomorphism $\rho_{r,q}$ induces the morphism
$$
p_{X,r,q}:J^r_0(\Bbb C^p,X)\to J^q_0(\Bbb C^p,X).
$$
In particular, we have the morphism $p_{X,q,0}:J^q_0(\Bbb C^p,X)\to X$.

The projective limit $J^\infty_0(\Bbb C^p,X)=\varprojlim J^q_0(\Bbb C^p,X)$
is called the {\it space of $\infty$-jets} at $0\in\Bbb C^p$ of morphisms
from $\Bbb C^p$ to $X$ or the {\it space of formal morphisms} from
$\Bbb C^p$ to $X$. 
By the definition of a projective limit we have natural
projections $p_{X,\infty,q}:J^\infty_0(\Bbb C^p,X)\to J^q_0(\Bbb C^p,X)$.
By definition, one can consider a point of $J^\infty_0(\Bbb C^p,X)$ either 
as a morphism  $\Bbb C[X]\to \Bbb C[[X_1,\dots,X_p]]$ or as a 
morphism from the completion $\tilde\Cal O_x$ of the local
ring $\Cal O_x$ for some $x\in X$ to $\Bbb C[[X_1,\dots,X_p]]$. 

In particular, for the above presentation of $X$ as a Zariski-closed subset 
of $\Bbb C^m$, each $h\in J^\infty_0(\Bbb C^p,X)$ is uniquely defined by 
a morphism $h':\Bbb C[W_1,\dots,W_m]\to\Bbb C[[X_1,\dots,X_p]]$ with 
$h'(\varPhi_l)=0$ for $l=1,\dots,r$. It is defined by 
$h'(W_i)=h_i\in\Bbb C[[X_1,\dots,X_p]]$, where 
$\varPhi_l(h_i)=0$ for each $l=1,\dots,r$.

Consider $X$ as a complex analytic set and let $h:\Bbb C^p,0\to X$ be
a germ of a holomorphic map at $0\in\Bbb C^p$. Denote by
$\frak F_{\Bbb C^p, 0}$ and $\frak F_{X,x}$ the rings of germs of
holomorphic functions on $\Bbb C^p$
at $0$ and on $X$ at $x$ respectively. We may identify the ring
$\frak F_{\Bbb C^p, 0}$ with a subring of the ring $\Bbb C[[X_1,\dots,X_p]]$.
Consider the morphism
$h^*:\frak F_{X,x}\to\frak F_{\Bbb C^p, 0}\subseteq\Bbb C[[X_1,\dots,X_p]]$
induced by $h$. The restriction of $h^*$ to $\Cal O_{X,x}$, which
is denoted by $j^\infty_0\,h$, belongs to $J^\infty_0(\Bbb C^p,X)$ and
is called the $\infty$-jet of $h$ at $0$. Put
$j^q_0\,h:=p_{X,\infty,q}(j^\infty_0\,h)$ and call $j^q_0\,h$ the $q$-jet of $h$
at $0$.

Denote by $x_1,\dots,x_p$ the standard coordinates in $\Bbb C^p$.
Let $A=(a_1,\dots,a_s)$ be a multi-index, $W$ a finite dimensional complex
vector space, and $F:\Bbb C^p,0\to W$ a germ of a holomorphic map, i.e\.
$F\in W\otimes\frak F_{\Bbb C^p,0}$.
We denote by $\p_A$ a linear operator on $W\otimes\frak F_{\Bbb C^p,0}$
which is equal to the tensor product of the identical operator on $W$ and
the operator $\frac{\p^{|A|}}{\p x_{a_1}\dots\p x_{a_s}}$ on
$\frak F_{\Bbb C^p,0}$.
In particular, we have $\p_\emptyset F=F$ and we write
$\p_a F$ instead of $\p_{(a)} F$.

For the above presentation of $X$ as a Zariski-closed subset in $\Bbb C^m$
the holomorphic germ $h:\Bbb C^p,0\to X$ can be given by a holomorphic
map $F$ from a neighborhood of $0\in\Bbb C^p$ to $\Bbb C^m$ such that
$\varPhi_l\o F=0$ for each $l=1,\dots,r$. Denote by $\frak A_p$ the set
of all multi-indices $A=(a_1,\dots,a_s)$. By definition,
the $\infty$-jet $j^\infty_0\,h$ is uniquely determined
by the indexed set $(\p_A F(0))_{A\in\frak A_p}$ of complex numbers.
The complex numbers $\p_A F(0)$ satisfy the equations
$\p_A(\varPhi_l\o F)(0)=0$ for $A\in\frak A_p$ and do not depend
on the choice of $F$. Similarly, the $q$-jet $j^q_0\,F$ is determined
by the indexed set $(\p_A F(0))_{A\in\frak A_{p,q}}$ of complex
numbers satisfying the equations $\p_A(\varPhi_l\o F)(0)=0$ for all
$A\in\frak A_{p,q}$. The above considerations show that for a smooth point
$x\in X$ our notion of jets coincide with the usual one.

Note that the jet spaces of holomorphic functions and of regular functions
on affine varieties coincide.

Later we denote by $\p_A$ also the linear operator on
$W\otimes\Bbb C[[X_1,\dots,X_p]]$ which is equal to the tensor product of
the identical operator on $W$ and the operator
$\frac{\p^{|A|}}{\p x_{a_1}\dots\p x_{a_s}}$ on $\Bbb C[[X_1,\dots,X_p]]$.

\subhead\nmb.{2.5} \endsubhead
Consider a $G$-module $V$, the spaces of $q$-jets $J^q_0(\Bbb C^p,V)$,
and $J^q_0(\Bbb C^p,Z)$, and the sets of formal morphisms
$J^\infty_0(\Bbb C^p,V)$ and $J^\infty_0(\Bbb C^p,Z)$.
The projection  $\pi:V\to Z$ induces the morphism
$J^q_0(\Bbb C^p,\pi):J^q_0(\Bbb C^p,V)\to J^q_0(\Bbb C^p,Z)$
and the map
$J^\infty_0(\Bbb C^p,\pi):J^\infty_0(\Bbb C^p,V)\to J^\infty_0(\Bbb C^p,Z)$.

The standard action of the group $G$ on $\Bbb C[V]$ induces an action of $G$
on $J^q_0(\Bbb C^p,V)$ by automorphisms of $J^q_0(\Bbb C^p,V)$ as an affine
variety and on $J^\infty_0(\Bbb C^p,V)$.
The inclusion $\Bbb C[V]^G\subseteq\Bbb C[V]$, the morphism
$J^q_0(\Bbb C^p,\pi)$, and the map $J^\infty_0(\Bbb C^p,\pi)$
induce the morphism
$$
\pi^q:J^q_0(\Bbb C^p,V)/G\to J^q_0(\Bbb C^p,Z)
$$
and the map
$$
\pi^\infty:J^\infty_0(\Bbb C^p,V)/G\to J^\infty_0(\Bbb C^p,Z).
$$
Denote by $\bar J^q_0(\Bbb C^p,Z)$ the Zariski closure of
$\pi^q(J^q_0(\Bbb C^p,V)/G)$ in $J^q_0(\Bbb C^p,Z)$.

\proclaim{Proposition} The morphism
$\pi^q:J^q_0(\Bbb C^p,V)/G\to J^q_0(\Bbb C^p,Z)$ induces a birational
morphism of $J^q_0(\Bbb C^p,V)/G$ to $\bar J^q_0(\Bbb C^p,Z)$.
\endproclaim
\demo{Proof} The group $G$ acts freely on the open subset
$p_{V,q,0}^{-1}(V_0)\subseteq J^q_0(\Bbb C^p,V)$. Since all points of $V_0$ and
$Z_0$ are smooth, for each $v\in V_0$ the morphism $J^q_0(\Bbb C^p,\pi)$
induces a bijective map of $p_{V,q,0}^{-1}(v)$ onto
$p_{Z,q,0}^{-1}(\pi(v))$.
Thus the morphism $\pi^q$ maps the Zariski-open
subset $p_{V,q,0}^{-1}(V_0)/G$ of $J^q_0(\Bbb C^p,V)/G$ onto the Zariski-open
subset $p_{Z,q,0}^{-1}(Z_0)$ of $J^q_0(\Bbb C^p,Z)$ bijectively and
this implies the statement of the proposition.
\qed\enddemo

\subhead{\nmb.{2.6}}\endsubhead
Evidently we have the following bijections
$$
J^q_0(\Bbb C^p,V)=\on{Hom}(\Bbb C[V],\frak J_p^q) 
=\on{Lin}(V^*,\frak J^q_p)=
V\otimes\frak J^q_p,
$$
where $V^*$ is the dual vector space for $V$, $\on{Hom}$ means 
the set of morphisms in the category of $\Bbb C$-algebras, and $\on{Lin}$ 
means the set of linear mappings. So each
$h\in J^q_0(\Bbb C^p,V)=V\otimes(\frak J_p^q)$
can be written uniquely as $h=\sum_{A\in\frak A_{p,q}}h_A\otimes X_A$,
where $h_A\in V$.
Similarly, $J^\infty_0(\Bbb C^p,V)=V\otimes\Bbb C[[X_1,\dots,X_p]]$ and, for
$h\in J^\infty_0(\Bbb C^p,V)$, we have $h=\sum_{A\in\frak A_p}h_A\otimes X_A$,
where $h_A\in V$.

\proclaim{Proposition} The space of jets $J^q_0(\Bbb C^p,V)$  is isomorphic to 
the affine space $V^{\frak A_{p,q}}$.
Moreover, there are isomorphisms of $G$-modules
$$
J^q_0(\Bbb C^p,V)=V^{\frak A_{p,q}},\quad\text{and}\quad
J^\infty_0(\Bbb C^p,V)=V^{\frak A_{p}}
$$
where the $G$-action on the products is the diagonal one.
\endproclaim
\demo{Proof} The first statement follows from the definition of the structure
of the affine variety on $J^q_0(\Bbb C^p,V)$. The maps
$J^q_0(\Bbb C^p,V)\ni h\mapsto (h_A)_{A\in \frak A_{p,q}}$ and
$J^\infty_0(\Bbb C^p,V)\ni h\mapsto (h_A)_{A\in\frak A_p}$ give the required
isomorphisms of $G$-modules.
\qed\enddemo
Note that for
$$
A=(\undersetbrace \text{$r_1$ times} \to {1,\dots,1},\dots,
\undersetbrace \text{$r_p$ times}\to {p,\dots,p}),
$$
where $r_1,\dots,r_p\ge 0$, and for a germ of a holomorphic map
$F:\Bbb C^p,0\to V$ we have
$h=j^\infty_0\,F=\sum_A\frac1{A!}\p_A\,F(0)\otimes X_A$,
where $A!=r_1!\dots r_p!$. 
A similar formula is true for the $q$-jet $j^q_0\,F$.

\subhead{\nmb.{2.7}. The problem of lifting}\endsubhead
We consider the following problem.
Let $f:\Bbb C^p\to Z$ be a rational morphism
which is regular on a Zariski-open subset $U$ of $\Bbb C^p$. A rational
morphism $F:\Bbb C^p\to V$ which is regular on $U$ is called a
{\it regular lift} of $f$ if $\pi\o F=f$.

$$
\xymatrix{
  &V\ar@{->}[1,0]^{\pi}  \\
\CC^p\ar@{->}[0,1]^{f} \ar@{->}[-1,1]^{F}  & Z \\
}
$$
Similarly, let $U$ be a connected classically open subset $U$ of
$\Bbb C^p$ and let $U\to Z$ be a holomorphic map, i.e\. a morphism in the
category of complex analytic sets. A holomorphic map $F:U\to V$ is called a
{\it holomorphic lift} of $f$ if $\pi\o F=f$. If $f:\Bbb C^p,x\to Z$ is a germ at
$x\in\Bbb C^p$ of a holomorphic map from $\Bbb C^p$ to $Z$, a germ 
$F:\Bbb C^p,x\to V$ at $x$ of a holomorphic map from $\Bbb C^p$ to $V$ is called a
{\it local lift} of $f$ if $\pi\o F=f$.

Let $f\in J^\infty_0(\Bbb C^p,Z)$ be a formal morphism. A formal morphism
$F\in J^\infty_0(\Bbb C^p,V)$ is called a {\it formal lift} of $f$ if
$J^\infty_0(\Bbb C^p,\pi)\o F=f$.

The aim of this paper is to find conditions for the existence of all
these lifts.

A germ $f:\Bbb C^p,x\to Z$ at $x\in\Bbb C^p$ of a holomorphic map from $\Bbb C^p$ to $Z$ 
is called {\it quasiregular} if $f^{-1}(Z_0)$ meets any neighborhood of $x$.

By \nmb!{2.3}, if a germ $f$ is quasiregular, there is a choice of invariant coordinates 
$y_i$ such that for the corresponding $\tilde\De$ the composition 
$\tilde\De\o f$ does not vanish identically. 

A formal morphism $f:\Bbb C[Z]\to\Bbb C[[X_1,\dots,X_p]]$ from $\Bbb C^p$ to 
$Z$ is called {\it quasiregular} if there is no stratum $S$ of $Z$ of codimension 
$\ge 1$ such that $f$ factors through a morphism 
$\Bbb C[\bar S]\to\Bbb C[[X_1,\dots,X_p]]$, where $\bar S$ is the closure of $S$. 

We claim that {\it for a quasiregular formal morphism $f$ and for any choice of 
invariant coordinates $y_i$, for the corresponding $\tilde\De$ we have 
$f(\tilde\De)\ne 0$}. 
By \nmb!{2.3} the set of zeros of all $\De$ obtained 
from any choice of invariant coordinates coincides with the Zariski-closed subset 
$Z_{>0}$ of $Z$. If our claim is wrong, by Hilbert's Nullstellensatz the formal 
morphism  $f$ vanishes on the ideal $I(Z_{>0})=\sqrt{I(Z_{>0})}$ 
of $\Bbb C[W_1,\dots,W_m]$  which defines $Z_{>0}$. 
Consider the standard presentation of  
$I(Z_{>0})$ as an intersection of a finite set of prime ideals of $\Bbb C[Z]$ 
corresponding to the decomposition of $Z_{>0}$ into irreducible components. 
Since $\Bbb C[[X_1,\dots,X_p]]$ is an integral domain, $f$ vanishes on at least 
one of these prime ideals. But each of these prime ideals defines a 
component of $Z_{>0}$ and such components are the closures of strata of $Z$ of 
codimension $\ge 1$. This contradicts our assumption.  

Let $K$ be a subgroup of $G$ and let $V^K$ be a subspace of $V$ of fixed points 
of the action of $K$ on $V$. Let $H=N_G(K)$ be the normalizer of $K$ in $G$ consisting
of all elements of $G$ preserving $V^K$, and let $W=H/K$ be the corresponding
quotient group acting naturally on $V^K$. By definition, $\pi(V^K)$ is the closure 
$\bar S$ of a stratum $S$ of $Z$.  

Consider the natural map 
$\ka: V^K/W\to\pi(V^K)=\bar S$. It is evidently bijective 
and regular since the natural map $V/W\to V/G$ is regular and 
$\pi(V^K)=\bar S$ is a Zariski-closed subset of $Z$ as  
the projection $\pi$ is a finite morphism. Then the morphism $\ka$ is birational.
Denote by $\bar S_{\text{nor}}$ the set of all normal points of $\bar S$, i.e\., points 
$x\in\bar S$ such that the local ring $\Cal O_x(\bar S)$ is integrally closed. 
It is known that 
$\bar S_{\text{nor}}$ is a Zariski-open subset of $\bar S$ and 
$S\subseteq \bar S_{\text{nor}}$ since $S$ is smooth. 
Since the affine variety $V^K/W$ is normal, 
by Zariski's main theorem, the restriction $\ka|_{S_n}$ of $\ka$ to 
$\bar S_{\text{nor}}$ induces 
an isomorphism between the algebraic varieties 
$\pi^{-1}(\bar S_{\text{nor}})$ and $\bar S_{\text{nor}}$. 
 
Assume that, for a holomorphic map $f$ as above, $f(U)$ is contained in
$Z_{>0}$. Then $f(U)$ is contained in the closure $\bar S$ of a
stratum $S$ of $Z$ of codimension $\ge 1$.
Namely, let $f(U)\subseteq Z_{>i-1}$ for maximal $i$. Then there exists
$x\in U$ such that $f(x)$ is a point of some stratum $S$ of
codimension $i$; otherwise $f(U)\subseteq Z_{>i}$. If a regular function
$h\in\Bbb C[Z]$ vanishes on $S$ then $h\o f$ vanishes on an
open neighborhood of $x$ in $U$ and thus on the whole of $U$. So
$f(U)\subseteq\bar S$ and there is a subgroup $K$ of $G$ distinct from $G$ 
such that $f(U)\subseteq\bar S=\pi(V^K)$. 

It is clear that if each morphism $f$ of the above type (regular, holomorphic, or 
formal) from $\Bbb C^p$ to $V^K/W$ has a lift $F$ (regular, holomorphic, local, or 
formal), then the composition of $f$ with the morphism $\ka$ has the corresponding 
lift to $V$ which is the composition of $F$ with the inclusion $V^K\to V$. 

Conversely, if $f:\Bbb C,x\to Z$ is a germ of a holomorphic map at $x\in\Bbb C^p$ such 
that $f(x)\in\bar S_{\text{nor}}$, there is a unique germ of a holomorphic map 
$f':\Bbb C^p,x\to V^K/W$ such that $\ka\o f'=f$. Similarly, let  
$f:\Bbb C[Z]\to\Bbb C[[X_1,\dots,X_p]]$ be a formal morphism from $\Bbb C^p$ to $Z$ 
which can be extended to the morphism $\tilde\Cal O_z(Z)\to\Bbb C[[X_1,\dots,X_p]]$ 
for some $z\in\bar S_{\text{nor}}$. There is a unique formal morphism 
$f':\Bbb C[V/W]\to\Bbb C[[X_1,\dots,X_p]]$ such that $J^\infty_0(\Bbb C^p,\ka)(f')=f$. 
In both cases the lifting problem for $\pi:V\to V/G$ reduces to 
the corresponding one for $\pi_{V^K}:V^K\to V^K/W$. 

Although these arguments give 
nothing if the above points $f(x)$ or $z$ do not belong $\bar S_{\text{nor}}$, we have 
the following theorem. 

\proclaim{Theorem} 
\roster
\item Let $f$ be a holomorphic map from a classically open subset 
$U$ of $\Bbb C^p$ to $Z$, $S$ a stratum of maximal codimension 
such that $f(U)\subseteq\bar S$, and $f:\Bbb C^p,x\to Z$ be a germ of $f$ 
at some $x\in U$. Let $K$ be a subgroup of 
$G$ such that $\pi(V^K)=\bar S$, and let $W=N_G(K)/K$.
If the germ $f:\Bbb C^p,x\to Z$ at $x\in U$ has a local lift to $V$, then the 
germ $f':\Bbb C^p,x\to V^K/W$ of the map defined by $f$ is a quasiregular germ 
of a holomorphic map and this germ has a local lift to $V^K$. 
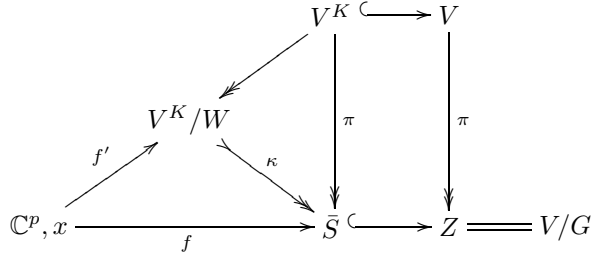
$$
\xymatrix{
  &  &  &{V^K\;}\ar@{->>}[2,0]^{\pi} \ar@{^(->}[0,1] \ar@{->>}[1,-1]  &V\ar@{->>}
[2,0]^{\pi}  &  \\
  &  &V^K/W\ar@{ >->>}[1,1]^{\ka}  &  &  &  \\
  &\CC ^p,x\ar@{->}[0,2]_{f} \ar@{->}[-1,1]^{f'}  &  &{\bar S\;}\ar@{^(->}[0,1] 
 &Z\ar@{={}}[0,1]  &V/G \\
}
$$
\item Let $f$ be a formal morphism from $\Bbb C^p$ 
to $Z$ and let $S$ be a stratum of the maximal codimension such that $f$ factors 
through a formal morphism $f'$ from $\Bbb C^p$ to $\bar S$. If $f$ has a formal lift 
to $V$, then the formal morphism $f'$ is quasiregular and has a formal lift to 
$V^K$ such that $\pi(V^K)=\bar S$. 
\item If $F_1$ and $F_2$ are holomorphic lifts of a holomorphic map $f:U\to Z$, then 
there is a $g\in G$ such that $F_2=g\o F_1$. The same is true for local lifts of 
germs of holomorphic maps, and for lifts of formal morphisms.
\endroster
\endproclaim
\demo{Proof} (1) Consider a local lift of $f$ which is a germ of  a holomorphic 
map $F:U'\to V$, where $U'\subseteq U$ is an open subset. By assumption,  
$F(U')\subseteq\pi^{-1}(\bar S)=\cap_{g\in G}gV^K$ and there is a point $x\in U'$ 
such that the stabilizer $G_{F(x)}=gKg^{-1}$. Then $F(U)\subseteq g(V^K)$ and         
$F'=g^{-1}\o F$ is a local lift of $f$ such that $F'(U')\subseteq V^K$. 
Then $\pi_{V^K}\o F'$ is a germ of a holomorphic map which by construction 
coincides with the germ of $f'$. By definition, the germ $f'$ is quasiregular 
and $F'$ is its local lift. 

(2) Let $F$ be a formal lift of $f$ to $V$. Let $I(\bar S)$ be the prime ideal of 
$\Bbb C[Z]$ defining $\bar S$. Consider the pullback $\pi^*(I(\bar S))$ of 
$I(\bar S)$. By the definition of $V^K$ we have $\pi^{-1}(\bar S)=\cup_{g\in G}gV^K$.
By definition of a formal lift the formal lift $F$ vanishes on $\pi^*(I(\bar S))$ and 
then, by Hilbert's Nullstellensatz, on the ideal 
$I(\cup_{g\in G}gV^K)=\sqrt{I(\cup_{g\in G}gV^K)}$ of $\Bbb C[V]$ defining 
a Zariski-closed subset $\cup_{g\in G}gV^K$ of $V$. 
Evidently the ideal $I(\cup_{g\in G}gV^K)$ 
equals the intersection of prime ideals $I(gV^K)$. Since $\Bbb C[[X_1,\dots,X_p]]$ is 
an integral domain, there is a $g\in G$ such that the formal morphism $F$ vanishes on 
$I(gV^K)$ and then $F\o g^{-1}$ is a formal lift of $f$ which factors through the formal 
morphism $F':\Bbb C[V^K]\to\Bbb C[[X_1,\dots,X_p]]$.
Thus the formal morphism $f$ factors through the 
formal morphism $f'=J^\infty_0(\Bbb C^p,\pi_{V^K})(F')$ from $\Bbb C^p$ to $V^K/W$,
$F'$ is a formal lift of $f'$, and, by assumption, the formal morphism $f'$ is 
quasiregular. 

(3) First assume that the germ of $f$ is quasiregular at $x\in U$. 

Let $F_1$ and $F_2$ be holomorphic lifts of $f$. By assumption, there is a point $y$ in 
a neighborhood of $x$ 
such that $F_1(y)$ and $F_2(y)$ are regular points of $V$. Since
$(\pi\o F_1)(y)=(\pi\o F_2)(y)$, there is a unique $g\in G$ such that
$F_2(y)=(g\o F_1)(y)$. As the projection $\pi$ is \'etale at $F_2(y)$,
the lift $F_2$ coincides with $g\o F_1$ in a neighborhood of $y$ and then on the
whole of $U$.

Let $K$ be the maximal subgroup of $G$ such that $f(U)\subseteq\pi(V^K)$, 
let $N_G(K)$ be the normalizer of $K$ in $G$, and $W=N_G(K)/K$.
By the proof of (1) the germ $f:\Bbb C^p,x\to Z$ for $x\in U$ can be considered 
as a quasiregular germ of a holomorphic map $f':\Bbb C^p,x\to\bar S=\pi(V^K)$ 
and there are $g_1,g_2\in G$ such that 
the germs $g_1\o F:\Bbb C^p,x\to V$  and $g_2\o F_2:\Bbb C^p,x\to V$ are local 
lifts of the 
above germ $f'$ to $V^K$. Then, for some $g\in H$ we have 
$g_2\o F_2=(gg_1)\o F_1$ in 
a neighborhood of $x$ and then in the whole of $U$. Thus we have 
$F_2=(g_2^{-1}gg_1)F_1$. 

For local lifts the proof is the same. For lifts of formal morphisms the
proof follows from \nmb!{5.4} below.
\qed\enddemo 

The above theorem shows that the problem of lifting for local and formal lifts 
is reduced in some sense to one for the quasiregular case. 

Namely, let
the conditions of theorem \nmb!{2.7}.1 be satisfied. Since the morphism 
$\ka$ is birational, for each basic invariant $\ta$ of 
$\Bbb C[V^K]$, the composition $\ka^{-1}\o\ta $ is a rational function 
on $V^K/W$ and, in general, the function $\ka^{-1}\o\ta\o f$ is a meromorphic 
function on $U$. First we have to check that this function is analytic 
near $x$. If $f(x)\in\bar S_{\text{nor}}$ this is always true, because $\ka^{-1}$ is 
an isomorphism near $f(x)$. 
Then $f$ has a local lift at $x$ iff the germ $f':\Bbb C^p,x\to V^K/W$  
has a local lift to $V^K$. 

The analogous statement for formal lifts is true whenever the conditions of 
theorem \nmb!{2.7}.2 are satisfied.     

\subhead{\nmb.{2.8}. An algebraic interpretation of the problem of lifting}
\endsubhead
The results of this section are not used in the rest of the paper.

The above geometric problem of lifting has the following algebraic
interpretation. For instance, suppose that
$f:\Bbb C^p\to Z$ is a regular morphism and $F$ is its regular lift. 
Consider the 
morphism $f^*:\Bbb C[Z]=\Bbb C[V]^G\to\Bbb C[\Bbb C^p]$ induced by $f$
and the morphism $F^*:\Bbb C[V]\to\Bbb C[\Bbb C^p]$ induced by $F$. 
Since, by definition, $\Bbb C[Z]\subseteq\Bbb C[V]$, 
the morphism $F^*$ is an extension of the morphism $f^*$ to 
$\Bbb C[V]$.

$$
\xymatrix{
  &V\ar@{->>}[1,0]  &  &\CC [V]\ar@{-->}[1,-1]^{F^*}  \\
\CC ^p\ar@{->}[0,1]^{f} \ar@{-->}[-1,1]^{F}  &Z &\CC [\CC ^p] &\CC [Z]\ar@
{_(->}[-1,0] \ar@{->}[0,-1]^{f^*}  \\
}
$$
Similarly, consider a germ of a holomorphic morphism $f:\Bbb C^p,0\to Z,O$,
where $O=\pi(0)$ and its local lift $F:\Bbb C^p,0\to V,0$. We have the morphisms
$f^*:\frak F_{Z,O}\to\frak F_{\Bbb C^p,0}$ and  
$F^*:\frak F_{V,0}\to\frak F_{\Bbb C^p,0}$ induced by $f$ and $F$ respectively.  
Since the projection $\pi$
induces the inclusion $\frak F_{Z,O}\subseteq\frak F_{V,0}$, the morphism $F^*$ 
is an extension of the morphism $f^*$ to $\frak F_{V,0}$. 

Finally, let $f:\Bbb C[Z]\to\Bbb C[[X_1,\dots,X_p]]$ be a formal morphism 
from $\Bbb C^p$ to $Z$ and let $F:\Bbb C[V]\to\Bbb C[[X_1,\dots,X_p]]$ be its 
formal lift. Since the projection $\pi$ induces the inclusion 
$\Bbb C[Z]\subseteq\Bbb C[V]$ the lift $F$ is an extension of 
the morphism $f$ to $\Bbb C[V]$.  

\subhead\nmb.{2.9} \endsubhead
Let $\ta$ be a homogeneous $G$-invariant polynomial of degree $d$ on $V$
and let $\ta^s$ be the corresponding  symmetric $d$-linear form on $V$.
For each germ $F:\Bbb C^p,0\to V$ of a holomorphic map and each system of
multi-indices $(A_1,\dots,A_d)$ we put
$$
T(A_1,\dots,A_d)(j^q_0\,F):=\ta^s\Bigl(\p_{A_1}F(0),\dots,\p_{A_d}F(0)\Bigr).
$$
 By \nmb!{2.4} $T(A_1,\dots,A_d)$ is a function on
$J^q_0(\Bbb C^p,V)$ for $q\ge |A_1|,\dots,|A_d|$.
 From \nmb!{2.2} and proposition \nmb!{2.6} it follows that 
the function 
$$
T(A_1,\dots,A_d):J^q_0(\Bbb C^p,V)=V^{\frak A_{p,q}}\to\Bbb C
$$ 
is regular, $G$-invariant, and equal to a polarization of
$\ta$ up to a nonzero factor. It is also symmetric in $A_1,\dots,A_d$.

By proposition \nmb!{2.5}, define a rational function
$\tilde T(A_1,\dots,A_d)$ on $\bar J^q_0(\Bbb C^p,Z)$ by the condition
$T(A_1,\dots,A_d)=\tilde T(A_1,\dots,A_d)\o J^q_0(\Bbb C^p,\pi)$.
By definition, we have
$$
\tilde T(\emptyset,\dots,\emptyset)\o\pi=T(\emptyset,\dots,\emptyset)=\ta.
\tag{\nmb:{1}}
$$

Now extend the $d$-linear form $\ta^s$ on $V$ to a $d$-linear form $\frak T$
on $J^\infty_0(\Bbb C^p,V)=V\otimes\Bbb C[[X_1,\dots,X_p]]$ with values in 
$\Bbb C[[X_1,\dots,X_p]]$ which is 
defined by the following condition. For $i=1,\dots,d$, $v_i\in V$,
and $F_i\in\Bbb C[[X_1,\dots,X_p]]$
$$
\frak T(v_1\otimes F_1,\dots,v_d\otimes F_d)
:=\ta^s(v_1,\dots,v_d)\;F_1\cdots F_d.
$$
For $h\in J^\infty_0(\Bbb C^p,V)=V\otimes\Bbb C[[X_1,\dots,X_p]]$ and
a system of multi-indices $A_1,\dots,A_d$ put
$$
\frak T(A_1,\dots,A_d)(h):=\frak T(\p_{A_1}h,\dots,\p_{A_d}h).
$$
By definition, the function 
$\frak T(A_1,\dots,A_d):J^\infty_0(\Bbb C^p,V)\to\Bbb C[[X_1,\dots,X_p]]$ 
is $G$-invariant and symmetric in $A_1,\dots,A_d$.

$$
\xymatrix{
 J_0^\oo (\CC ^p,V) \ar@{={}}[0,1]  &
V\otimes\CC[[X_1,\dots,X_p]] \ar@{->}[0,2]^{\frak T(A_1,\dots,A_d)}&  &\CC [[X_1,\dots ,X_p]] \\
  &J_0^q(\CC ^p,V)\ar@{->}[0,2]^{T(A_1,\dots,A_d)} \ar@{->}[1,0]^
{J_0^q(\CC^p,\pi)}  &  &\CC  \\
  &\bar J_0^q(\CC ^p,Z)\ar@{->}[-1,2]_{\qquad\tilde T(A_1,\dots,A_d)}  &  &  \\
}
$$
\head\totoc\nmb0{3}. The functions $\tilde T(A_1,\dots,A_d)$ and
$\tilde\frak T(A_1,\dots,A_d)$ \endhead

\subhead{\nmb.{3.1}. The $q$-jet of the identity map on $V$ in
invariant coordinates}\endsubhead

Let $v:V\to V$ be the identity map. Let $v_0\in V_0$ be a regular point of $V$ 
and let $y_i$ be invariant
coordinates in a neighborhood $U$ of $v_0$ in $V$ introduced in \nmb!{2.3}.
Then in $U$ the map $v$ is defined by a holomorphic function $v(y_i)$ with values in $V$.

Let $I=(i_1,\dots,i_s)$ be a (unordered) multi-index with $i_1,\dots,i_s\in\{1,\dots,n\}$. 
In particular, for $s=0$ we put $I:=\emptyset$. Then the $q$-jet of the identity map $v$ 
at each point $x\in U$ is defined by the set of partial derivatives 
$\p_Iv=\frac{\partial^sv}{\partial y_{i_1}\dots\partial y_{i_s}}$ for $|I|\le q$ at $x$. 

Let $e_a$ be a basis of $V$, let $u_a$ be the corresponding coordinates, and let 
$J=\on{det}(\frac{\partial y_i}{\partial u_j})$ be the jacobian. 

\proclaim{Lemma} Let $I=(i_1,\dots,i_s)$, where $s>0$, be a multi-index.
Then $\widetilde{\p_Iv}:=J^{2s-1}\p_Iv$ is a regular map from $U$ to $V$. 
\endproclaim
\demo{Proof}
We prove this lemma by induction with respect to $s$.
We use that
$$
\frac{\partial u_a}{\partial y_i}=\frac1{J}J^a_i,
$$
where $J^a_i\in\Bbb C[V]$ is the algebraic complement of the entry
$\p y_i/\p u_a$ in the Jacobi matrix.

For $s=1$, $J\p_iv=\sum_{a=1}^nJ^a_ie_a$ is a regular map from $U$ to $V$.

Let $I=(i_1,\dots,i_s)$ with $\widetilde{\p_{I}v}$ being regular.
Then for $I'=(i_1,\dots,i_{s+1})$ we have
$$
\p _{I'}v=\frac{\p}{\p y_{i_{s+1}}}\,\Bigl(\frac{\widetilde{\p_Iv}}{J^{2s-1}}\Bigr)=
\frac1{J^{2s+1}}\sum_{a=1}^nJ^a_{i_{s+1}}\Bigl(J\frac{\p(\widetilde{\p_Iv})}
{\p u_a}-
(2s-1)\frac{\p J}{\p u_a}\widetilde{\p_Iv}\Bigr),
$$
where $\sum_{a=1}^nJ^a_{i_{s+1}}\Bigl(J\frac{\p(\widetilde{\p_Iv})}{\p u_a}-
(2s-1)\frac{\p J}{\p u_a}\widetilde{\p_Iv}\Bigr)$
is a regular map from $U$ to $V$.
\qed\enddemo

\subhead{\nmb.{3.2}}\endsubhead
Let $h=j^q_0\,F\in J^q_0(\Bbb C^p,V)$, where $F:\Bbb C^p,0\to V$ is a
germ of a holomorphic map such that $F(0)\in V_0$.
Put $F_i:=y_i\o F$, where $y_i$ are the invariant coordinates on $V$.
We need to express the $q$-jet $j_0^qF$ in terms of the $q$-jet of the identity 
map $v$, i.e\., we have to find the explicit formula for each
$h_A=\p_A F(0)$
with $A\in\frak A_{p,q}$ in terms of $\p_B F_i$ and $\p_Iv$
with $|B|,|I|\le |A|$. We can extract this formula from
the following expression (see the classical Faa di Bruno formula for $p=1$
and \cit!{18} or \cit!{19} for arbitrary $p$).
$$
d^q(v\o F)=q!\sum_{k=1}^q\frac{1}{k!}\sum_{i_1,\dots,i_k=1}^n
\Bigl(\frac{\partial^kv}{\partial y_{i_1}\dots\partial y_{i_k}}\o F\Bigr)
\sum_{q_1+\dots+q_k=q}\frac{d^{q_1}F_{i_1}}{q_1!}\dots
\frac{d^{q_k}F_{i_k}}{q_k!}.\tag{\nmb:{1}}
$$
where $y_i$ are arbitrary local coordinates in $V$. Note that the formula 
\thetag{\nmb|{1}} is true whenever $F:\Bbb C[V]\to\Bbb C[[X_1,\dots,X_p]]$ 
is a formal morphism from $\Bbb C^p$ to $V$ and $F_i=F(y_i)$. 

The formula \thetag{\nmb|{1}} implies the following

\proclaim{Lemma} For each multi-index $A=(a_1,\dots,a_s)\ne\emptyset$,
where $a_1,\dots,a_s\in\{1,\dots,p\}$,
there is a well defined function
$$
\varPsi_A:(X,Y)
\mapsto \sum_{1\le |I|\le |A|}a_{A,I}(Y)\,X_I,
$$
where $X=(X_I)_{1\le |I|\leq |A|}$ and where the coefficients $a_{A,I}$
are polynomials in $Y=(y_{i,B})_{1\leq i\leq n,1\le |B|\leq |A|}$,
such that for each germ of a holomorphic map $F:\Bbb C^p,0\to V,v$ with
$v$ regular and for the local coordinates $y_i$ from above we have
$$
\p_A F=\varPsi_A\Bigl((\p_I v\o F)_{I},(\p_B F_i)_{i,B}\Bigr).
$$
\endproclaim

For example,
$$
\p_aF=\sum_{i=1}^n(\p_iv\o F)\;\p_a F_i,
$$
$$
\p_{(a_1,a_2)}F=\sum_{i,j=1}^n(\p_{(i,j)}v\o F)\, \p_{a_1}F_i\,\p_{a_2} F_j+
\sum_{i=1}^n(\p_iv\o F)\,\p_{(a_1,a_2)}F_i,
$$
and so on.

\subhead\nmb.{3.3} \endsubhead
We consider $T(A_1,\dots,A_d)$ which is a regular function
on $J^q_0(\Bbb C^p,V)$,
and $\tilde T(A_1,\dots,A_d)$ which is a rational function
on $\bar J^q_0(\Bbb C^p,Z)$, both defined in \nmb!{2.9}.

Let $z_i$ be the regular function on $Z$ for $i=1,\dots,n$,
used in \nmb!{2.3} for the construction
of the invariant coordinates on $V$. Let $A_1,\dots,A_{d'}\ne\emptyset$ and
$A_{d'+1}=\dots=A_d=\emptyset$. Put $M:=2(|A_1|+\dots+|A_d|)-d'$.
By lemma \nmb!{3.1}, for any system of multi-indices $I_1,\dots,I_d$ such that
$1\le |I_1|\le |A_1|,\dots,1\le |I_{d'}|\le |A_{d'}|$ and
$I_{d'+1}=\dots=I_d=\emptyset$, the expression
$\Delta_M\cdot T\o (\p_{I_1}v,\dots,\p_{I_d}v)$ is a $G$-invariant and regular
function on $V$. Thus there is a unique rational function
$\tilde T(I_1,\dots,I_d)$ on $Z$ such that $\tilde T(I_1,\dots,I_d)\o\pi=
T\o (\p_{I_1}v,\dots,\p_{I_d}v)$ and $\tilde\De_M\cdot \tilde T(I_1,\dots,I_d)$ is a
regular function on $Z$.

Let $q$ be the maximal order of the multi-indices $A_1,\dots,A_d$. For
$k=1,\dots,d'$ we may consider the $a_{A_k,I_k}$ of \nmb!{3.2}
as polynomials in
$Y=(y_{i,B})_{1\le i\le n, 1\le |B|\leq q}$.
Put
$$
a_{A_1,\dots,A_{d'},I_1,\dots,I_{d'}}(Y):=
a_{A_1,I_1}(Y)\cdots a_{A_{d'},I_{d'}}(Y).
$$

\proclaim{Theorem} Let $A_1,\dots,A_{d'}\ne\emptyset$,
$A_{d'+1},\dots,A_d=\emptyset$ and $A_1,\dots,A_d\in\frak A_{p,q}$. Then
\roster
\item Then
the following formula
defines a
rational function
$\tilde T(A_1,\dots,A_d)$ on $\bar J^q_0(\Bbb C^p,Z)$:
$$
\tilde T(A_1,\dots,A_d):=
\sum_{\Sb 1\le |I_1|\le |A_1|,\\
\dotsc \\1\le |I_{d'}|\le |A_{d'}|\endSb}
a_{A_1,\dots,A_{d'},I_1,\dots,I_{d'}}\Bigl((\p_Bz_i)_{i,B}\Bigr) \;
\tilde T(I_1,\dots,I_{d'},\emptyset,\dots,\emptyset);
$$
\item $\tilde T(A_1,\dots,A_d)\o J^q_0(\Bbb C^p,\pi)=T(A_1,\dots,A_d)$;
\item $\tilde\De_M\cdot \tilde T(A_1,\dots,A_d)$ is a regular function
on $J^q_0(\Bbb C^p,Z)$.
\endroster
\endproclaim
$$
\xymatrix{
J_0^q(\CC^p,V)\ar@{->}[0,3]^{\quad T(A_1,\dots,A_d)} \ar@{->}[1,0]_
{J_0^q(\CC^p,\pi)}  &  &  & \CC \\
\bar J_0^q(\CC^p,Z)\ar@{-->}[-1,3]_{\qquad\tilde T(A_1,\dots,A_d)}  &  &  & \\
}       
$$
\demo{Proof} (1) By proposition \nmb!{2.5}, it suffices to check that the
condition $\tilde T(A_1,\dots,A_d)\o\pi^q=T(A_1,\dots,A_d)$ is satisfied
for the above expression of $\tilde T(A_1,\dots,A_d)$. By lemma \nmb!{3.2},
this follows from
$$\multline
T(A_1,\dots,A_d)(h)=\\
\sum_{\Sb 1\le |I_1|\le |A_1|,\\ \dotsc \\1\le |I_{d'}|\le |A_{d'}|\endSb}
\Bigl(a_{A_1,\dots,A_{d'},I_1,\dots,I_{d'}}\bigl((\p_B F_i)_{i,B}\bigr)\;
T(\p_{I_1}v,\dots,\p_{I_{d'}}v,v,\dots,v)\o F\Bigr)(0),
\endmultline
$$
where $h=j^q_0\,F\in J^q_0(\Bbb C^p,V)$.

(2) This statement follows from proposition \nmb!{2.5}.

(3) This statement follows from (1) and lemma \nmb!{3.1}.
\qed\enddemo

\subhead\nmb.{3.4} \endsubhead
Let $F\in V\otimes\Bbb C[[X_1,\dots,X_p]]$ be a formal morphism from
$\Bbb C^p$ to $V$ and $F_i=y_i(F)$.
lemma \nmb!{3.2} implies the following

\proclaim{Lemma} For each multi-index $A$ such that $|A|\ge 1$
and the invariant coordinates $y_i$ on $V$ we have
$$
J^{2|A|-1}(F)\;\p_A F=
\varPsi_A\biggl(\Bigl((J^{2(|A|-|I|)}\widetilde{\p_I v})(F)\Bigr)_{I},
\Bigl(\p_B F_i\Bigr)_{i,B}\biggr).
$$
\endproclaim

\subhead{\nmb.{3.5}. The function $\tilde\frak T(A_1,\dots,A_d)$} \endsubhead
Consider the presentation of $Z$ as an irreducible  Zariski-closed subset
of $\Bbb C^m$ defined in \nmb!{2.1}.
Denote by $I(Z)$ the prime ideal of the ring of polynomials
$\Bbb C[W_1,\dots,W_m]$ defining $Z\subseteq\Bbb C^m$. By \nmb!{2.4}, each
formal morphism $f\in J^\infty_0(\Bbb C^p,Z)$ is defined by the equations
$f(W_j)=f_j$ for $(j=1,\dots,m)$, where $f_j\in\Bbb C[[X_1,\dots,X_p]]$ 
and $\varPhi(f_j)=0$ for each $\varPhi\in I(Z)$.

Let $\psi$ be a regular function on $Z$ which is the restriction to
$Z$ of a polynomial $\Psi\in\Bbb C[W_1,\dots,W_m]$.
For each $f\in J^\infty_0(\Bbb C^p,Z)$  put $\psi(f):=\Psi(f_j)$. 
By definition we have $\psi(f)=f(\psi)$, where $f$ is considered as 
a morphism $\Bbb C[Z]\to\Bbb C[[X_1,\dots,X_p]]$. Then  
$\psi(f)$ defines a unique function 
$J^\infty_0(\Bbb C^p,Z)\to\Bbb C[[X_1,\dots,X_p]]$ 
which is independent of the choice of the polynomial $\Psi$.

Similarly, consider a rational function $\psi$ on $Z$ such that
$\psi=\frac{\psi_1}{\psi_2}$, where $\psi_1$ and $\psi_2$ are regular
functions on $Z$ and put $\psi(f):=\frac{\psi_1(f)}{\psi_2(f)}$
whenever $\psi_2(f)\ne 0$.
It is clear that $\psi(f)$ is a function on $J^\infty_0(\Bbb C^p,Z)$
with values in the field $\Bbb C((X_1,\dots,X_p))$ of fractions of the ring
$\Bbb C[[X_1,\dots,X_p]]$ which is independent of the choice
of the presentation  $\psi=\frac{\psi_1}{\psi_2}$.

Let $z_i$ be the regular functions on $Z$ used in \nmb!{2.3} for the construction
of the invariant coordinates $y_i$ on $V$.
Let $f\in J^\infty_0(\Bbb C^p,Z)$ be a quasiregular formal morphism from
$\Bbb C^p$ to $Z$ such that $\tilde\De(f)\ne 0$.
For $A_1,\dots,A_{d'}\ne\emptyset$,
$A_{d'+1},\dots,A_d=\emptyset$ and $M=2(|A_1|+\dots+|A_d|)-d'$ put
$$\multline
\tilde\frak T(A_1,\dots,A_{d'},\emptyset,\dots,\emptyset)(f):=\\
=\sum_{\Sb 1\le |I_1|\le |A_1|,\\ \dotsc \\1\le |I_{d'}|\le |A_{d'}|\endSb}
\Bigl(a_{A_1,\dots,A_{d'},I_1,\dots,I_{d'}}\bigl((\p_Bz_i)_{i,B}\bigr)\cdot
\tilde T(I_1,\dots,I_{d'},\emptyset,\dots,\emptyset)\Bigr)(f).
\endmultline$$
By definition, $\tilde{\frak T}(A_1,\dots,A_d)$ is a function with values in
the field $\Bbb C((X_1,\dots,X_p))$ on the set $\bar J^\infty_0(\Bbb
C^p,Z)$ consisting of all $f\in J^\infty_0(\Bbb C^p,Z)$ 
such that $\tilde\De(f)\ne 0$.

\proclaim{Theorem} Let $A_1,\dots,A_{d'}\ne\emptyset$ and
$A_{d'+1},\dots,A_d=\emptyset$ and $M=2(|A_1|+\dots+|A_d|)-d'$.
The function $\tilde{\frak T}(A_1,\dots,A_d)$  
satisfies the following
conditions:
\roster
\item $\tilde\frak T(A_1,\dots,A_d)\o J^\infty_0(\Bbb C^p,\pi)=
\frak T(A_1,\dots,A_d)$, where $\frak T$ is from \nmb!{2.9}.
\item $\tilde\De_M\cdot\tilde\frak T(A_1,\dots,A_d)$, where $\tilde\De_M$
is regarded as a function on $J^\infty_0(\Bbb C^p,Z)$, is a function
on $J^\infty_0(\Bbb C^p,Z)$ with values in $\Bbb C[[X_1,\dots,X_p]]$.
\endroster
\endproclaim
$$
\xymatrix{
 J_0^\oo(\CC ^p,V)\ar@{->}[0,3]^{\frak T(A_1,\dots,A_d)\quad} \ar@{->}[1,0]_
{J_0^\oo(\CC^p,\pi)}  &  & &\CC[[X_1,\dots,X_p]]  \\
\bar J_0^\oo(\CC ^p,Z)\ar@{-->}[-1,3]_{\qquad\tilde \frak T(A_1,\dots,A_d)}  
&  & & \\}       
$$
\demo{Proof} The proof follows from the definition of the function
$\tilde\frak T(A_1,\dots,A_d)$ and lemma \nmb!{3.4}.
\qed\enddemo

\head\totoc\nmb0{4}. The conditions of local and global lifting\endhead

\subhead{\nmb.{4.1}}\endsubhead
First we consider local lifts at regular points.

\proclaim{Proposition} Let $f:\Bbb C^p,x\to Z,z$ be a germ at
$x\in\Bbb C^p$ of a holomorphic map with $z$ regular. Then for each
$v\in\pi^{-1}(z)$ there is a unique local holomorphic lift
$F:\Bbb C^p,x\to V,v$ of $f$.
\endproclaim
\demo{Proof} By proposition \nmb!{2.1}, the map
$\pi$ is \'etale on $V_0$. Thus for each point $v\in\pi^{-1}(z)$
there is a unique local holomorphic lift $F:\Bbb C^p,x\to V_0,v$ of $f$.
\qed\enddemo

\subhead{\nmb.{4.2}. Lifts of quasiregular holomorphic germs}\endsubhead
Let $X$ be an affine variety and let $f$ be either a rational
morphism from $\Bbb C^p$ to $X$ or a holomorphic map defined on a
classically open connected subset $U\subseteq \Bbb C^p$ to $X$.
Consider the morphism $j^qf$ from $\Bbb C^p$ or from $U$ to
$J^q_0(\Bbb C^p,X)$,
which for $x\in U$ is given by $j^qf(x)=j_0^qf(\quad+x)$.
The morphism $j^qf$ is rational and is regular wherever
$f$ is regular; or holomorphic if $f$ is holomorphic.

Let $\si:V\to\si(V)\subseteq\Bbb C^m$ be the morphism defined by the system of
basic generators $\si_1,\dots,\si_m$.
Recall that $\si(V)$ and $Z=V/G$ are isomorphic as affine varieties and,
for this presentation of $Z$, $\si$ equals $\pi$.

Denote by $w_1,\dots,w_m$ the standard coordinates in $\Bbb C^m$ and let
$I(Z)$ be the prime ideal of the ring $\Bbb C[W_1,\dots,W_m]$ defining 
$Z$. Consider $\Bbb C[W_1,\dots,W_m]$ as a graded ring with a
grading defined by $\deg W_j=\deg\si_j$ for $j=1,\dots,m$.
Then $I(Z)$ is a homogeneous ideal.

For $\ta=\si_j$ and a system $A_1,\dots,A_{d_j}$ of multi-indices
denote by $\tilde S_j(A_1,\dots,A_{d_j})$ the rational function
$\tilde T(A_1,\dots,A_{d_j})$ from \nmb!{2.9} on
$\bar J^q_0(\Bbb C^p,Z)$.

Recall that by \nmb!{2.7} for a quasiregular germ $f:\Bbb C^p,0\to Z$ of 
a holomorphic map there is a choice of invariant coordinates such that 
for the corresponding function $\tilde\De$ we have $\tilde\De\o f\ne 0$.  
 
\proclaim{Theorem} Consider a quasiregular germ $f:\Bbb C^p,0\to Z$ of
a holomorphic map described by $w_j\o f=f_j$ for $j=1,\dots,m$.
Assume that, for some choice of the invariant coordinates such that
$\tilde\De\o f\ne 0$, $q$ is the minimal order
of nonzero terms of the Taylor expansion of $\tilde\De\o f$ at $0$.

Then the lift $F$ of $f$ at $0$ exists if and only if for 
$j=1,\dots, m$ and for each system of multi-indices
$A_1,\dots,A_{d_j}\in\frak A_{p,q}$
the functions $\tilde S_j(A_1,\dots,A_{d_j})\o j^qf$ have holomorphic
extensions to a neighborhood of $0$.
\endproclaim
$$
\xymatrix{
 &  &  &J_0^q(\CC ^p,V)\ar@{->}[0,3]^{\quad S(A_1,\dots,A_d)} \ar@{->}[1,0]_
{J_0^q(\CC^p,\pi)\!\!}  &  & & \CC  \\
\CC ^p,0\ar@{->}[0,3]^{j^qf} \ar@{->}[-1,3]^{j^qF}  & &  &\bar J_0^q(\CC ^p,Z)
\ar@{->}[-1,3]_{\qquad\tilde S(A_1,\dots,A_d)}  &  &  &\\
}       
$$
\demo{Proof} Let $F$ be a lift of $f$. Then by the definition of the
function $\tilde S_j(A_1,\dots,A_{d_j})$ for each $q\ge 0$ we have
$$
\tilde S_j(A_1,\dots,A_{d_j})\o j^qf=S_j(A_1,\dots,A_{d_j})\o j^q F:
\Bbb C^p,0\to \Bbb C,
$$
where the right hand side defines a holomorphic germ.

Conversely, let the assumptions of the theorem be satisfied.
Let us now use a representative $f:U\to Z$ of the germ, where $U$ is a
connected open neighborhood of $0$. Let $q$ be
the minimal order of nonzero terms of the Taylor expansion of
$\tilde\De\o f$ at $0$.
For each $j=1,\dots, m$ consider the function
$$
f^q_j(x,t):=\sum_{A_1,\dots,A_{d_j}\in\frak A_{p,q}}
\Bigl(\tilde S(A_1,\dots,A_{d_j})\o j^qf\Bigr)(x)\;
t_{A_1}\cdots t_{A_{d_j}},
\tag{\nmb:{1}}
$$
where $x\in U$ and $t=(t_A)_{A\in \frak A_{p,q}}\in\Bbb C^{\frak A_{p,q}}$.
By assumption, the function $f^q_j$ is a polynomial in $t$ whose coefficients
are holomorphic near 0.
By definition, the map
$f^q=(f^q_j)_{j=1,\dots,m}:\Bbb C^p\times\Bbb C^{\frak A_{p,q}}\to \Bbb C^m$
is holomorphic near 0.

Since $Z_{>0}$ is a Zariski-closed subset of $Z$ of
codimension $\ge 1$, the inverse image $f^{-1}(Z_{>0})$ is a complex analytic
subset of $U$ of codimension $\ge 1$ and $f^{-1}(Z_0)$ is a dense open
subset of $U$.

Let, for $y\in U$, $f(y)$ be a regular point and let $F_y$ be a local lift of
$f$ defined in a neighborhood $U_y$ of $y$, which exists by
proposition \nmb!{4.1}. For each $q$ consider the holomorphic map
$F^q_y:U_y\times\Bbb C^{\frak A_{p,q}}\to V$ given by:
$$
F^q_y(x,t):=\sum_{A\in\frak A_{p,q}}\p_AF_y(x)\;t_A.
\tag{\nmb:{2}}
$$
By theorem \nmb!{3.3}, we have
$$
\align
(\si_j\o F^q_y)(x,t)&=\sum_{A_1,\dots,A_{d_j}\in\frak A_{p,q}}
S_j\Bigl(\p_{A_1}F_y(x),\dots,\p_{A_{d_j}}F_y(x)\Bigr)\;
t_{A_1}\cdots t_{A_{d_j}}\\ &=
\sum_{A_1,\dots,A_{d_j}\in\frak A_{p,q}}
\Bigl(S_j(A_1,\dots,A_{d_j})\o j^qF_y\Bigr)(x)\;t_{A_1}\cdots t_{A_{d_j}}
=f^q_j(x,t).
\endalign$$
Therefore for each polynomial $\varPhi\in I(Z)$ we have
$\varPhi\o f^q=0$ on $U_y\x\Bbb C^{\frak A_{p,q}}$ and thus also on $U\x\Bbb C^{\frak A_{p,q}}$.
So $f^q$ is a holomorphic map from $U\x\Bbb C^{\frak A_{p,q}}$ to $Z$ and
$F^q_y$ is a lift of $f^q$.

For each germ of a holomorphic function $\varphi\in\frak F_{\Bbb C^p,x}$,
denote by $\on{Tay}_x^q\varphi$ the sum of terms of the Taylor
expansion at $x$ of $\varphi$ of orders $\le q$. For each germ
$\varphi=(\varphi_j)_j$ of a holomorphic map $\Bbb C^p,x\to\Bbb C^m$,
put $\on{Tay}^q_x\varphi:=(\on{Tay}^q_x\varphi_j)_j$.

By assumption, there is a multi-index $A\in\frak A_{p,q}$ such that
$\p_A(\tilde\De\o f)(0)=\p_A(\tilde\De\o\on{Tay}_0^qf)(0)\neq 0$.
This implies that there is a point $x_0=(x_{0,1},\dots,x_{0,p})\in\Bbb C^p$
such that $(\tilde\De\o \on{Tay}_0^qf)(x_0)\neq 0$.

For
$$
A=(\undersetbrace \text{$r_1$ times} \to {1,\dots,1},\dots,
\undersetbrace \text{$r_p$ times}\to {p,\dots,p}),
$$
put
$$
t_A(x):= \frac1{r_1!\dots r_p!}\;(x_1)^{r_1}\dots (x_p)^{r_p},\quad
t(x):=(t_A(x))_{A\in \frak A_{p,q}}
$$
where $x=(x_1,\dots,x_p)\in\Bbb C^p$.

By definition, we have $F^q_y(y,t(x-y))=\on{Tay}_y^qF_y(x)$ and then
$f^q_j(y,t(x-y))=(\si_j\o\on{Tay}_y^qF_y)(x)$.
On the other hand, since $\si_j$ is homogeneous, for a fixed $y$ we have
$\on{Tay}^q_yf_j=\on{Tay}^q_y(\si_j\o F_y)=
\on{Tay}_y^q(\si_j\o\on{Tay}^q_yF_y)$.
Thus, we have
$$
\on{Tay}^q_yf^q_j(y,t(x-y))=\on{Tay}^q_yf_j(x).\tag{\nmb:{3}}
$$
By assumption, the function $\tilde S_j(A_1,\dots,A_{d_j})\o j^qf$ has
a holomorphic extension to a neighborhood of $0$ and we may suppose that
the point $y$ belongs to this neighborhood. Letting $y\to 0$ in \thetag{\nmb|{3}}
we get $\on{Tay}^q_0f^q_j(0,t(x))=\on{Tay}^q_0f_j(x)$. Then we have
$(\tilde\De\o f^q)(0,t(x))=\tilde\De(\on{Tay}^q_0f)(x)\neq 0$ and, for the
point $x_0\in\Bbb C^p$ chosen above, we have
$(\tilde\De\o f^q)(0,t(x_0))\neq 0$, i.e\., $f^q(0,t(x_0))$ is a
regular point of $Z$.

Now we will construct a local lift of $f$. Consider a local holomorphic
lift $F^q$ of $f^q$ near $(0,t(x_0))$ in $U\times \Bbb C^{\frak A_{p,q}}$ which
exists by proposition \nmb!{4.1}. We can choose $y$ near $0$ so that
$f(y)\in Z_0$, so there exists a local holomorphic lift $F_y$ of $f$ near $y$,
and still $(\tilde\De\o f^q)(y,t(x_0))\ne 0$.
Consider the map $F^q_y$ defined by formula \thetag{\nmb|{2}}. Both
$F_y^q$ and $F^q$ are local lifts of $f^q$ at $(y,t(x_0))$.
By theorem \nmb!{2.7}, there exists $g\in G$ such that $F^q_y=gF^q$ near
$(y,t(x_0))$.

Since $F_y^q(x,t)$ is linear in $t\in \Bbb C^{\frak A_{p,q}}$, also $F^q(x,t)$ is
linear in $t$ and thus is defined for all $t$. Put $t_1:=(t_{1,A})_A$, where
$t_{1,\emptyset}=1$ and $t_{1,A}=0$ for $A\ne \emptyset$. Then
near $0\in\Bbb C^p$
we have by \nmb!{2.9.1},
$$
\si_j(F^q(x,t_1))=f^q_j(x,t_1)=
\Bigl(\tilde S_j(\emptyset,\dots,\emptyset)\o j^qf\Bigr)(x)=f_j(x),
$$
i.e\., $F^q(\quad,t_1)$ is a local lift of $f$ at $0$.
\qed\enddemo

\remark{Remark}
Consider the grading of the ring $\Bbb C[Z]=\Bbb C[V]^G$ induced by the
natural grading of the polynomial ring $\Bbb C[V]$ and denote by $r$
the order of the homogeneous function $\tilde\De$.
Let $f:\Bbb C^p,0\to Z$ be a germ of a holomorphic map satisfying
for some positive integer $q$ and for each $j=1,\dots,m$
the following conditions:
\roster
\item The function $\tilde S_j(A_1,\dots,A_{d_j})\o j^qf$ has a holomorphic
extension to a neighborhood of $0$ for each system of multi-indices
$A_1,\dots,A_{d_j}\in\frak A_{p,q}$ such that,
$\Bigl(\tilde S_j(A_1,\dots,A_{d_j})\o j^{q-1}f\Bigr)(0)=0$
for all $A_1,\dots,A_{d_j}\in\frak A_{p,q-1}$;
\item $\on{Tay}^{rq}_0(\tilde\De(f))\ne 0$.
\endroster
Then $f$ has a local lift at $0$.

Actually, since $\on{Tay}^{rq}_0(\tilde\De(f))=
\on{Tay}^{rq}_0(\tilde\De(\on{Tay}^q_0f))$, the proof of theorem \nmb!{4.2}
is valid for this $q$.
\endremark

\subhead{\nmb.{4.3}. The conditions for lifting of first order} \endsubhead
Next we use the notion of rational tensor fields on affine varieties
(see \cit!{7}).

For each $0\le s\le d$ consider the rational symmetric
tensor field $\ta_s$ of type $\binom{0}{s}$ on $Z$ defined as follows:
$$
\ta_s:=
\sum_{i_1,\dots,i_s=1}^n\tilde T\Bigl((i_1),\dots,(i_s),\emptyset,\dots,\emptyset\Bigr)
\,dz_{i_1}\otimes\dots\otimes dz_{i_s}
$$
where $\tilde T\Bigl((i_1),\dots,(i_s),\emptyset,\dots,\emptyset\Bigr)$ is a partial
case of the function $\tilde T(I_1,\dots,I_d)$ defined in
\nmb!{3.3}.

Consider the pull back $\pi^*\ta_s$ of $\ta_s$. By definition, we have
$$
\pi^*\ta_s=\sum_{i_1,\dots,i_s=1}^n
T(v_{i_1},\dots,v_{i_s},v,\dots,v)\,dy_{i_1}\otimes\dots\otimes dy_{i_s}=
T(dv,\dots,dv,v,\dots,v)
$$
and then $\ta_s=\pi_*T(dv,\dots,dv,v,\dots,v)$, where $dv$ occurs precisely
$s$ times. Since the projection
$\pi$ is \'etale on $V_0$ and by the above formula for $\pi^*\ta_s$ the
tensor field $\pi^*\ta_s$ is independent of the choice of the invariant
coordinates $y_i$, the tensor field $\ta_s$ is independent of the choice
of the invariant coordinates as well.
Note that $\ta_1=\frac1{d}d\ta$ induces a regular differential $1$-form on $Z$.

By lemma \nmb!{3.2}, for each germ $f:\Bbb C^p,0\to Z$ of a holomorphic map
we have
$$
f^*\ta_s=\sum_{a_1,\dots,a_s}^p
\Bigl(\tilde T\bigl((a_1),\dots,(a_s),\emptyset,\dots,\emptyset\bigr)\o j^1f\Bigr)\cdot
dx_{a_1}\otimes\dots\otimes dx_{a_s}.
$$
Denote by $\si_{j,s}$, where $(1\le s\le d_j)$, the tensor field $\ta_s$ for
$\ta=\si_j$.
Then the conditions for lifting of theorem \nmb!{4.2} for each
$j=1,\dots,m$ and $A_1,\dots,A_{d_j}\in\frak A_{p,1}$ are equivalent to the
following statement: For each $1\le s_j\le d_j$
the pull back $f^*\si_{j,s_j}$ is a holomorphic germ of a symmetric
covariant tensor field on $\Bbb C^p$.
We call these conditions the conditions of the first order. For $s=1$
these conditions are satisfied automatically.

Note that the conditions of the remark at the end of \nmb!{4.2} for $q=1$
use the conditions of first order only. For example, it suffices to
consider only these conditions if we need to have lifts which are
linear maps from $\Bbb C^p$ to $V$.

\subhead{\nmb.{4.4}. Global holomorphic lifts} \endsubhead
The following theorem shows that the problem of global holomorphic lifting
can be described topologically.

\proclaim{Theorem} Let  $U\subseteq\Bbb C^p$ be a classically open
connected subset of $\Bbb C^p$ and let $f:U\to Z$ be a holomorphic
map such that $f^{-1}(Z_0)\neq\emptyset$. Then a holomorphic lift
$F:U\to V$ exists
iff the image of the fundamental group $\pi_1(f^{-1}(Z_0))$ under $f$ is
contained in the image of the fundamental group $\pi_1(V_0)$ under the
projection  $\pi$.
\endproclaim

\demo{Proof} Since by proposition \nmb!{4.1}, the local holomorphic lift of
$f$ exists for each $x\in f^{-1}(Z_0)$, the condition of the theorem is
equivalent to the existence of a holomorphic lift for the restriction of $f$
to $f^{-1}(Z_0)$. Actually, let $F$ be such a lift. Since $f^{-1}(Z_0)$ is
an open dense subset of $U$ and $\pi$ is a finite morphism, the lift $F$
is bounded on bounded subsets of $U\cap f\i(Z_0)$.
Then by Riemann's extension theorem $F$ has a holomorphic
extension to $U$ which is a holomorphic lift of $f$.
\qed\enddemo

\subhead\nmb.{4.5}\endsubhead
We indicate that the problem of the existence of a global regular lift
reduces to the one for a holomorphic lift.

\proclaim{Theorem} Let $U\subseteq\Bbb C^p$ be a Zariski-open subset of
$\Bbb C^p$ and let $f:\Bbb C^p\to Z$ be a rational morphism which is
regular in $U$ and such that $f^{-1}(Z_0)\ne\emptyset$.

If a global holomorphic lift of $f$ on $U$ exists then it is regular.
\endproclaim

\demo{Proof} The proof follows from lemma 5.1.1 of \cit!{7}.
\qed\enddemo

\subhead{\nmb.{4.6}. Global regular lifts}\endsubhead
Now we indicate conditions for the existence of a global regular lift.

\proclaim{Theorem} Let $f:\Bbb C^p\to Z$ be a regular morphism
such that $f(\Bbb C^p)\cap Z_0\ne\emptyset$. Then $f$ has a regular lift
iff there is an integer $q>0$ such that, for each $j=1,\dots m$ and
each multi-index $A=(a_1,\dots,a_q)$,
the mapping $\tilde S_{j}(A,\dots,A)\o j^qf$ is constant.
\endproclaim

\demo{Proof} Let $u_i$ be linear coordinates in $V$ and let 
$F=(F_1,\dots,F_n)$ be
the expression of a regular lift of $f$ in these coordinates.
Suppose $q$ is the maximal degree of the polynomials $F_i$. Then
for each $A=(a_1,\dots,a_q)$ the map $S_j(A,\dots,A)\o j^qF$ is constant.
Theorem \nmb!{3.3} implies that $\tilde S_{j}(A,\dots,A)\o j^qf$ is
constant as well.

Let the condition of the theorem be satisfied. Let $x\in\Bbb C^p$ be
a point such that $f(x)\in Z_0$. Then there is a local
lift $F$ of $f$ at $x$. By theorem \nmb!{4.2}, the condition of the theorem
implies that $\si_j(\p_AF)$ is constant for each $A=(a_1,\dots,a_q)$
and each $j$. But this means that $\p_AF$ is constant also and, therefore,
$F$ is a polynomial map of degree $\le q$ in a neighborhood of $x$. Thus $F$
has a polynomial extension to the whole of $\Bbb C^p$ and this extension
is a lift of $f$.
\qed\enddemo

\subhead\nmb.{4.7} \endsubhead
We consider a special case when the existence of a local lift implies the
existence of a global lift.

\proclaim{Theorem} Let $f=(f_j):\Bbb C^p\to\Bbb C^m$ be
a regular morphism from $\Bbb C^p$ to $Z$ such that
each function $f_j$ is homogeneous of degree $rd_j$ for some positive
integer $r$. Then a global regular lift of $f$ exists iff
$f$ has a local holomorphic lift at $0\in\Bbb C^p$.
\endproclaim
\demo{Proof} Consider the action of the group $\Bbb C^*$ on $Z$ induced
by the action of the homothety group on $V$. This action induces a homotopy
equivalence between the open subset $f^{-1}(Z_0)$ and
the open subset $f^{-1}(Z_0)\cap B$, where $B$ is an open ball in
$\Bbb C^p$ centered at $0$. Then the statement of the theorem follows from
theorems \nmb!{4.4} and \nmb!{4.5}.
\qed\enddemo

\head\totoc\nmb0{5}. Formal lifts\endhead
In this section we find the conditions for lifts of quasiregular formal
morphisms from $\Bbb C^p$ to $Z$. Note first that proposition \nmb!{4.1}
about the existence of lifts at regular points also holds for formal
morphisms.

Let $\si_1,\dots,\si_m$ be the basic generators of
$\Bbb C[V]^G$, $\deg\si_j=d_j$ for $j=1,\dots,m$,
and $\si:V\to\si(V)\subseteq\Bbb C^m$ the corresponding morphism.
Consider some invariant coordinates $y_i$ on $V$ and the corresponding
function $\tilde\De$ on $Z$. Recall that we consider $Z$ as a
Zariski-closed subset $\si(V)$ of $\Bbb C^m$ defined by the ideal
$I(Z)$ of the ring $\Bbb C[W_1,\dots,W_m]$ and, for this presentation of $Z$, the
projection $\pi$ equals the map $\si:V\to\si(V)\subseteq\Bbb C^m$.

\subhead\nmb.{5.1} The functions $P_q(\ta)$ and $\tilde P_q(\ta)$\endsubhead
For a homogeneous $G$-invariant polynomial $\ta$ on $V$,
consider the function 
$P_q(\ta):J^\infty_0(\Bbb C^p,V)\to\Bbb C[(t_A)_A]\otimes\Bbb C[[X_1,\dots,X_p]]$ 
and the function $\tilde P_q(\ta)$ on the set of quasiregular formal
morphisms $f\in J^\infty_0(\Bbb C^p,Z)$ such that $\tilde\De(f)\neq 0$
with values in $\Bbb C[(t_A)_A]\otimes\Bbb C((X_1,\dots,X_p))$, 
where $(t_A)_A=(t_A)_{A\in\frak A_{p,q}}$ and 
$\Bbb C[(t_A)_A]$ is the ring of polynomials 
in $(t_A)_A$ with complex coefficients, defined as follows:
$$\gather
P_q(\ta)(F):=\sum_{A_1,\dots,A_d\in \frak A_{p,q}}\frak T(A_1,\dots,A_d)(F)\;
t_{A_1}\dots t_{A_d},\\
\tilde P_q(\ta)(f):=\sum_{A_1,\dots,A_d\in \frak A_{p,q}}\tilde \frak T(A_1,\dots,A_d)(f)\;
t_{A_1}\dots t_{A_d},
\endgather$$
where $F\in J^\infty_0(\Bbb C^p,V)$, $f\in J^\infty_0(\Bbb C^p,Z)$, and
where $\frak T(A_1,\dots,A_d)$ and $\tilde \frak T(A_1,\dots,A_d)$ 
are the functions 
defined in \nmb!{2.9} and \nmb!{3.5}. 

The following lemma follows from the definitions of $P_q(\ta)$,
$\tilde P_q(\ta)$, $\tilde\frak T(A_1,\dots,A_d)$, and theorem \nmb!{3.5}.

\proclaim{Lemma}
\roster
\item  We have
$$
P_q(\ta)=\tilde P_q(\ta)\o J^\infty_0(\Bbb C^p,\pi): 
J^\infty_0(\Bbb C^p,V) \to \Bbb C[(t_A)_A]\otimes\Bbb C[[X_1,\dots,X_p]].
$$
$$
\xymatrix{
 J_0^\oo(\CC ^p,V)\ar@{->}[0,2]^{P_q(\ta)\qquad\quad} \ar@{->}[1,0]_
{J_0^\oo(\CC^p,\pi)}  &  &\CC[(t_A)_A]\otimes \CC[[X_1,\dots,X_p]]  \\
\bar J_0^\oo(\CC ^p,Z)\ar@{->}[-1,2]_{\quad\tilde P_q(\ta)}  &  &  \\
}       
$$
\item If $\ta_1,\ta_2\in\Bbb C[V]^G$ are homogeneous polynomials of the same
degree, then $\tilde P_q(\ta_1+\ta_2)=\tilde P_q(\ta_1)+\tilde P_q(\ta_2)$;
\item Let $\ta_1,\ta_2\in\Bbb C[V]^G$ be homogeneous polynomials. Then we
have  $\tilde P_q(\ta_1\ta_2)=\tilde P_q(\ta_1)\tilde P_q(\ta_2)$.
\item Let $f$ be a polynomial in the graded variables $T_1,\dots,T_r$
of degrees $d_1,\dots,d_r$ which is homogeneous with respect to
this grading, and let
$\ta_1,\dots,\ta_r\in\Bbb C[V]^G$ be homogeneous polynomials of degrees
$d_1,\dots,d_r$. Then we have
$$\gather
P_q\Bigl(f(\ta_1,\dots,\ta_r)\Bigr)=f\Bigl(P_q(\ta_1),\dots,P_q(\ta_r)\Bigr), \\
\tilde P_q\Bigl(f(\ta_1,\dots,\ta_r)\Bigr)=f\Bigl(\tilde P_q(\ta_1),\dots,\tilde P_q(\ta_r)\Bigr).
\endgather$$
\endroster
\endproclaim

\subhead\nmb.{5.2}  \endsubhead
For $\si=(\si_1,\dots,\si_m)$ and a formal morphism 
$F:\Bbb C[V]\to\Bbb C[[X_1,\dots,X_p]]$ from $\Bbb C^p$ to $V$ put 
$\si(F):=(F(\si_1),\dots,F(\si_m))$.

\proclaim{Lemma} Let $F:\Bbb C[V]\to\Bbb C[[X_1,\dots,X_p]]$ be a formal 
morphism from $\Bbb C^p$ to $V$. If $\si(F)=0$, then $F$ vanishes on the
set of all regular functions on $V$ with zero constant terms; or, $F=0$ as
an element of $V\otimes C[[X_1,\dots,X_p]]$.
\endproclaim

\demo{Proof} Let $(e_i)$ be a basis of $V$ and $u_i$ the corresponding 
coordinates. It is sufficient to prove that $F(u_i)=0$. Since the group $G$ 
is finite the ring $\Bbb C[V]$ is integral over its subalgebra $\Bbb C[V]^G$.
Then for each $i=1,\dots,n$ there is a polynomial 
$p(x)=x^N+\sum_{j=1}^Na_{N-j}x^{N-j}$, whose coefficients $a_{N-i}$ belong 
to $\Bbb C[V]^G$, such that $p(u_i)=0$. Consider the natural grading of the 
ring $\Bbb C[V]$. Since $\deg ((u_i)^N)=N$ we may assume that $\deg a_{N-j}=j$.
This implies that the coefficients $a_{N-j}$ as polynomials in $\si_j$ have no 
constant terms. Then we have 
$$
0=F(p(u_i))=F(u_i)^N+\sum_{j=1}^nF(a_{N-j})F(u_i)^{N-j}.
$$
Since $\si(F)=0$, this equation implies $F(u_i)^N=0$ and therefore $F(u_i)=0$. 
\qed\enddemo

\subhead{\nmb.{5.3}. The conditions for formal lifts}\endsubhead
For $F\in J^\infty_0(\Bbb C^p,V)$ consider $P_q(\ta)(F)$ as a polynomial
in $(t_A)_A$ with coefficients in $\Bbb C[[X_1,\dots,X_p]]$.
Denote by $P_q(\ta)(F)_0(t)$ the polynomial in $(t_A)_A$ which is obtained
by the evaluation of the coefficients of the polynomial
$P_q(\ta)(F)$ at $X=(X_1,\dots,X_p)=0$.
Similarly, for $f\in J^\infty_0(\Bbb C^p,Z)$ consider $\tilde P_q(\ta)(F)$ as
a polynomial in $(t_A)_A$ with coefficients in $\Bbb C((X_1,\dots,X_p))$
and denote by $\tilde P_q(\ta)(f)_0(t)$ the polynomial in $t$ which is
obtained by the evaluation of the coefficients of the polynomial
$\tilde P_q(\ta)(f)$ at  $X=(X_1,\dots,X_p)=0$ whenever their values at 
$X=0$ are defined.

For
$$
A=(\undersetbrace \text{$r_1$ times} \to {1,\dots,1},\dots,
\undersetbrace \text{$r_p$ times}\to {p,\dots,p}),
$$
put
$$
t_A(X):= \frac1{r_1!\dots r_p!}\;(X_1)^{r_1}\dots (X_p)^{r_p},
\quad t(X):=(t_A(X))_{A\in\frak A_{p,q}}.
$$
For a formal power series $\varphi\in\Bbb C[[X_1,\dots,X_p]]$, denote by
$\on{Tay}^q\varphi$ the sum of the terms of $\ph$ of orders $\le q$.
Denote by $\tilde{\frak S}_j(A_1,\dots,A_{d_j})$
the function $\tilde{\frak T}(A_1,\dots,A_{d_j})$ for $\ta=\si_j$.

Recall that by \nmb!{2.7} for a quasiregular formal morphism 
$f\in J^\infty_0(\Bbb C^p,Z)$ 
there is a choice of invariant coordinates such that for the corresponding function 
$\tilde\De$ we have $f(\tilde\De)\ne 0$ which here we write also as 
$\tilde\De(f)\ne 0$.  

\proclaim{Theorem} Let $f\in J^\infty_0(\Bbb C^p,Z)$ be a quasiregular
formal morphism given by the equations $f(w_j)=f_j\in\Bbb C[[X_1,\dots,X_p]]$ for
$j=1,\dots,m$, where $w_j$ are the standard coordinate functions on $\Bbb C^m\supseteq Z$.
Let $y_i$ be invariant coordinates on $V$ such that for
the corresponding function $\tilde\De$ we have $\tilde\De(f)\ne 0$.
Assume $q$ is the minimal order of nonzero terms of $\tilde\De(f)$.

Then a formal lift $F$ of $f$ exists iff
for $j=1,\dots, m$ and for each system of multi-indices
$A_1,\dots,A_{d_j}\in \frak A_{p,q}$ we have
$\tilde{\frak S}_j(A_1,\dots,A_{d_j})(f)\in\Bbb C[[X_1,\dots,X_p]]$ and
$\on{Tay}^qf_j=\tilde P_q(\si_j)(f)_0(t(X))$.
\endproclaim

\demo{Proof} Let $F$ be a formal lift of $f$. Then, by theorem \nmb!{3.5},
we have
$$
\tilde{\frak S}_j(A_1,\dots,A_{d_j})(f)=\frak S_j(A_1,\dots,A_{d_j})(F)
\in\Bbb C[[X_1,\dots,X_p]].
$$
Moreover, by lemma \nmb!{5.1}, we have
$$
\on{Tay}^q f_j=\on{Tay}^q(\si_j(F))=P_q(\si_j)(F)_0(t(X))
=\tilde P_q(\si_j)(f)_0(t(X)).
$$

Conversely, let the assumptions of the theorem be satisfied and let $q$ be
the minimal order of nonzero terms of $\tilde\De(f)$. By assumption,
$\on{Tay}^q(\tilde\De(f))\ne 0$. Then there is a point
$x_0=(x_{0,1},\dots,x_{0,p})\in\Bbb C^p$ such that
$\on{Tay}^q(\tilde\De(f))(x_0)\neq 0$.

For each $j=1,\dots, m$ consider the function
$$
f^q_j((t_A)_A):=\tilde P_q(\si_j)(f)=\sum_{A_1,\dots,A_{d_j}\in \frak A_{p,q}}
\tilde\frak S_j(A_1,\dots,A_{d_j})(f)\,
t_{A_1}\dots t_{A_{d_j}}.
$$
We may consider $f^q=(f^q_j)$ as a formal morphism
$\Bbb C^{\frak A_{p,q}}\x\Bbb C^p,(t(x_0)),0)\to\Bbb C^m$, i.e\. as a
morphism $\tilde\Cal O_{(\Bbb C^m,f^q(t(x_0),0)}\to
\tilde\Cal O_{\Bbb C^{\frak A_{p,q}}\x\Bbb C^p,(t(x_0),0)}$.
We prove that $f^q=(f^q_j)$ is a formal morphism
$\Bbb C^{\frak A_{p,q}}\x\Bbb C^p,(t(x_0),0)\to Z$ by the following arguments.
Let $\varPhi\in I(Z)$ be a homogeneous polynomial. Then $\varPhi\o\si=0$
and, by lemma \nmb!{5.1}, we have
$\varPhi(f^q_j)=\varPhi(\tilde P_q(\si_j)(f))=
\tilde P_q(\varPhi\o\si_j)(f)=0$.

By assumption, we have $\on{Tay}^qf_j=\tilde P_q(\si_j)(f)_0((t_A(X))_A)
=f^q_j((t_A(X))_A,0)$.
Since the polynomial $\De$ is homogeneous, we have
$$
\on{Tay}^q\Bigl(\tilde\De(f^q(t(X),0))\Bigr)(x_0)=
\on{Tay}^q(\tilde\De(\on{Tay}^qf))(x_0)=
\on{Tay}^q(\tilde\De(f))(x_0)\ne 0.
$$
Thus, the formal morphism $f^q=(f^q_j)_{j=1,\dots,m}$ has a formal lift
$$
F^q:\Bbb C^{\frak A_{p,q}}\x\Bbb C^p,((t(x_0),0)\to V,
$$
which can be written as follows: $F^q=\sum_AF^q_A\,t_A$, where
$A$ is a multi-index and $F^q_A\in V\otimes\Bbb C[[X_1,\dots,X_p]]$.

Since $F^q$ is a formal lift of $f^q$, for each $j=1,\dots,m$ we have
$$\gather
\si_j(F^q)=\sum_{A_1,\dots,A_{d_j}}\frak S_j(F^q_{A_1},\dots,F^q_{A_{d_j}})\;
t_{A_1}\dots t_{A_{d_j}}\\=
\sum_{A_1,\dots,A_{d_j}\in \frak A_{p,q}}\tilde\frak S_j(A_1,\dots,A_{d_j})
(f^q)\;(t_{A_1}+t_{A_{1}}(x_0))\dots(t_{A_{d_j}}+t_{A_{d_j}}(x_0)),
\endgather$$
This implies that $\frak S_j(F^q_{A_1},\dots,F^q_{A_{d_j}})=0$ whenever for some
$k=1,\dots, d_j$ we have $|A_k|>q$. In particular, for a multi-index
$A$ such that $|A|>q$ and for each $j=1,\dots, m$ we have
$\si_j(F^q_A)=0$. By lemma \nmb!{5.2},
we have $F^q_A=0$ and, therefore, the formal lift $F^q$ is a polynomial
in $(t_A)_A$ with coefficients in $V\otimes\Bbb C[[X_1,\dots,X_p]]$. Put $t_1:=(t_1^A)_A$ where
$t_1^\emptyset=1$, and $t_1^A=0$ for $A\ne\emptyset$. Denote by $F^q(t_1)$ the
value of $F^q$ as a polynomial in $t$ at $t=t_1$. Then we have
$$
\si_j(F^q(t_1))=f^q_j(t_1)=\tilde\frak S_j(\emptyset,\dots,\emptyset)(f)=f_j,
$$
i.e\. $F^q(t_1)$ is a formal lift of $f$.
\qed\enddemo

\subhead{\nmb.{5.4}}\endsubhead
Theorem \nmb!{5.3} implies the following
\proclaim{Corollary} The map
$\pi^\infty:J^\infty_0(\Bbb C^p,V)/G\to J^\infty_0(\Bbb C^p,Z)$ is
injective.
\endproclaim
\demo{Proof} Let $f\in J^\infty_0(\Bbb C^p,Z)$ be a formal
morphism which has a lift to $V$.

First assume that the morphism $f$ is quasiregular.
Consider a formal morphism
$f^q=(f^q_j)_j$ from $\Bbb C^{\frak A_{p,q}}\x\Bbb C^p,(t(x_0),0)$ to $Z$
constructed for $f$ in the proof of theorem \nmb!{5.3} and one of its lifts
$F^q:\Bbb C^{\frak A_{p,q}}\x\Bbb C^p,(t(x_0),0)\to V$. Since $F^q(t(x_0),0)$ is
a regular point of $V$, the lift $F^q$ is defined up to the action of some
$g\in G$. On the other hand, for each lift $F$ of $f$,
$F^q=\sum_A\frac1{A!}\p_AF\;t_A$ is a lift of $f^q$. This implies
that the lift $F$ of $f$ is defined up to the action of some $g\in G$.

For an arbitrary formal morphism $f\in J^\infty_0(\Bbb C^p,Z)$, there is
a subgroup $K$ of $G$ such that we can consider $f$ as a quasiregular
formal morphism to $V^K/(N_G(K)/K)$.
Then one can prove our statement using the same arguments as in the proof
of theorem \nmb!{2.7}.
\qed\enddemo

\head\totoc\nmb0{6}. Complex reflection groups\endhead

In this section we consider the case when $G$ is a finite group generated
by complex reflections.

\subhead{\nmb.{6.1}. Relative invariants}\endsubhead
Recall some standard facts about finite complex reflection groups
(see \cit!{15}, \cit!{16}, and \cit!{17}).

First note that in this case the basic generators $\si_1,\dots,\si_m$
are algebraically independent, $m=n=\dim V$, and then $Z=\Bbb C^n$. Therefore
we have a unique (up to permutation) choice of invariant coordinates $y_i$ .

Suppose $\frak H$ is the set of reflection hyperplanes of $G$ and, for
each $H\in\frak H$, denote by $e_H$ the order of the cyclic subgroup $G_H$
of $G$ fixing $H$ pointwise, by $s_H$ a generator of $G_H$,
and by $l_H$ a linear functional with the kernel $H$.
Then the Jacobian $J$ equals
$\prod_{H\in\frak H}l_H^{e_H-1}$ up to some nonzero constant factor $c$ and
one can take $\prod_{H\in\frak H}l_H^{e_H}$ for $\De$.  We may choose
the functionals $l_H$ such that $J=\prod_{H\in\frak H}l_H^{e_H-1}$.

Denote by $E_G$ the set of the orders $e_H$ of the cyclic subgroups $G_H$
fixing pointwise the reflection hyperplanes and, for each $e\in E_G$, set
$\frak H_e:=\{H\in\frak H|e_H=e\}$. Then $\De_e:=\prod_{H\in\frak H_e}l_H^e$
is $G$-invariant polynomial. Denote by $\tilde\De_e$ the regular
function on $Z$ such that $\tilde\De_e\o\pi=\De_e$.

Consider the standard action of the group $G$ on $\Bbb C[V]$ given by $g\cdot f=f\o g^{-1}$
for every $f\in\Bbb C[V]$ and $g\in G$.

Let $\chi:G\to\Bbb C\setminus 0$ be a character of $G$. A $\chi$-{\it relative
invariant} is a polynomial $f\in\Bbb C[V]$ such that $g\cdot f=\chi(g)f$ for each $g\in G$.

\proclaim{Theorem} Let $G$ be a finite group generated by complex
reflections and $\chi$ a character of $G$ such that
$\chi(s_H)=\on{det}^{-k_s}(s_H)$ for some $k_s$. Then the polynomial
$f_{\chi}:=\prod_{H\in\frak H}l_{H}^{k_s}$ is a $\chi$-relative invariant
and the space of $\chi$-relative invariants equals the space
$\Bbb C[V]^Gf_{\chi}$.
\endproclaim
In particular, the polynomial $J_e=\prod_{H\in\frak H_e}l_H^{e-1}$ is
a $\chi_e$-relative invariant for a character $\chi_e$ such that
$\chi_e(s_H)=\on{det}^{1-e}(s_H)$, when $H\in\frak H_e$, and $\chi_e(s_H)=1$,
when $H\notin\frak H_e$. Then the jacobian $J$ is a $\chi$-relative invariant
for a character $\chi$ such that $\chi(s_H)=\on{det}^{1-e_H}(s_H)$.

\subhead{\nmb.{6.2}. The functions $\tilde T(A_1,\dots,A_d)$ and
$\tilde\frak T(A_1,\dots,A_d)$ for complex reflection groups}\endsubhead
For a complex reflection group $G$ one can improve the statement (3) of
theorem \nmb!{3.3} and the statement (2) of theorem \nmb!{3.5} as follows.

\proclaim{Theorem} Let $\ta$ be a homogeneous $G$-invariant polynomial
of degree $d$ on $V$ and let $A_1,\dots,A_d$ be a system of multi-indices
such that $A_1,\dots,A_{d'}\neq\emptyset$ and $A_{d'+1}=\dots=A_d=\emptyset$.
Put $M:=2(|A_1|+\dots+|A_d|)-d'$ and let, for each $e\in E_G$, both $\mu_e$ and
$\nu_e$ be positive integers such that $M=\mu_ee-\nu_e$, where
$0\le\nu_e\le e$.

Then
\roster
\item $\prod_{e\in E_G}\tilde\De_e^{M-\mu_e}\tilde T(A_1,\dots,A_d)$ is a
regular function on $Z$;
\item Regard $\tilde\De_e$ as a function on $J^\infty_0(\Bbb C^p,Z)$.
Then $\prod_{e\in E_G}\tilde\De_e^{M-\mu_e}\;\tilde\frak T(A_1,\dots,A_d)$
is a function on $J^\infty_0(\Bbb C^p,Z)$ with values in
$\Bbb C[[X_1,\dots,X_p]]$.
\endroster
\endproclaim
\demo{Proof} (1) Consider the $G$-invariant rational function
$T(\p_{I_1}v,\dots,\p_{I_d}v)$ on $V$, where $I_1,\dots,I_d$ be a system of
multi-indices such that $1\le |I_1|\le |A_1|,\dots,1\le |I_{d'}|\le |A_{d'}|$ and
$I_{d'+1}=\dots=I_d=\emptyset$.
Since $\De_e=\prod_{H\in\frak H_e}l_H^e$ is a $G$-invariant polynomial,
by lemma \nmb!{3.1}, $J^M\;T(\p_{I_1}v,\dots,\p_{I_d}v)$ is a $\chi$-relative
invariant for the character $\chi$ where
$\chi(s_H)=\on{det}(s_H)^{-\nu_e}$ for $H\in\frak H_e$. By theorem
\nmb!{6.1} we have
$J^M\;T(\p_{I_1}v,\dots,\p_{I_d}v)\in\Bbb C[V]^Gf_{\chi}$, where
$f_{\chi}=\prod_{e\in E_G}\prod_{H\in\frak H_e}l_H^{\nu_e}$. 
On the other hand, we have $J^M=\prod_{e\in E_G}\De_e^{M-\mu_e}f_{\chi}$. 
Therefore
$\prod_{e\in E_G}\De_e^{M-\mu_e}\;T(I_1,\dots,I_d)$
is a $G$-invariant polynomial on $V$ and (1) follows from the definition
of the function $\tilde T(I_1,\dots,I_d)$ and theorem \nmb!{3.3} (1).

One can prove (2) similarly using the definition of the function
$\tilde\frak T(A_1,\dots,A_d)$ and lemma \nmb!{3.4}.
\qed\enddemo

To apply the above results on lifting we need to find explicit expressions
for the functions $\tilde S_j(A_1,\dots,A_{d_j})$ and
$\tilde\frak S_j(A_1,\dots,A_{d_j})$ for particular $G$-modules $V$. Although
this problem is purely technical the computations are rather complicated.
Therefore we consider only simple examples.

\subhead{\nmb.{6.3}. First Example}\endsubhead
Consider the simplest case when $V=\Bbb C$ and $G$ is the complex reflection
group generated by a generalized reflection $z\mapsto\on{exp}(\frac{2\pi i}{n})z$ for some
fixed $n\geq 2$. There is one basic invariant $\si:z\mapsto z^n$ and $Z=\Bbb C$.
Consider, for example, the lifting problem for $f\in\frak F_{\Bbb C^p,0}$,
i.e\. the problem of solving the equation $z^n=f$ in the ring
$\frak F_{\Bbb C^p,0}$. Note that
in this case $f$ is either quasiregular, or equals $0$,
and it suffices to assume that $f$ is quasiregular.

Let $y=z^n$ be the invariant coordinate on $V$. It is clear that the symmetric
$n$-linear form $S$ on $V$ corresponding to $\si$ equals $z_1\cdots z_n$.
Consider the system of multi-indices $A_1,\dots,A_n$ where
$A_1=(a_1^1,\dots,a_{q_1}^1),\dots,A_r=(a_1^r,\dots,a_{q_r}^r)$ and
$A_{r+1}=\dots=A_n=\emptyset$. Put
$$
f^r_{a_1^1\dots a^1_{q_1},\dots,a_1^r\dots a^r_{q_r}}:=
S(\p_{A_1}z,\dots,\p_{A_n}z)=z^{n-r}\p_{A_1}z\cdots\p_{A_r}z.
$$
By the general procedure we need to express
$f^r_{a_1^1\dots a^1_{q_1},\dots,a_1^r\dots a^r_{q_r}}$ via $y=z^n$ and its
partial derivatives $\p_Ay$. We do this by the following recurrence relations.
$$
\aligned
f^1_a&=\frac1n\p_ay,\\
y^{r-1}f^r_{a_1^1\dots a^1_{q_1},\dots,a_1^r\dots a^r_{q_r}}&=
f^1_{a_1^1\dots a^1_{q_1}}\cdots f^1_{a_1^r\dots a^r_{q_r}},\\
f^1_{a_1\dots a_q}&=\p_{a_q} f^1_{a_1\dots a_{q-1}}-
(n-1)f^2_{a_q,a_1\dots a_{q-1}},
\endaligned
\tag{\nmb:{1}}
$$
which give the required formulas by induction with respect to
$q=\operatorname{max}\{q_1,\dots,q_r\}$.

Thus, for the above system of multi-indices $A_1,\dots,A_n$ we have
$$
\tilde S(A_1,\dots,A_n)\o j^qf=
f^r_{a_1^1\dots a^1_{q_1},\dots,a_1^r\dots a^r_{q_r}},
$$
for $f\in J^q_0(\Bbb C,Z)$, and
$$
\tilde\frak S(A_1,\dots,A_n)(f)=
f^r_{a_1^1\dots a^1_{q_1},\dots,a_1^r\dots a^r_{q_r}}(f),
$$
for $f\in J^\infty_0(\Bbb C^p,Z)$.

There is the following stronger form of theorem \nmb!{4.2} for the case
under consideration.

\proclaim{Theorem} Let $f:\Bbb C^p,0\to Z$ be a germ of a holomorphic map.
Then $f$ has a local lift at $0$ iff either $f=0$, or there is a system of
indices $(a_1,\dots,a_r)$ such that the functions
$f^1_{a_1\dots a_r}\o j^rf$ and
$f^n_{a_1\dots a_r,\dots, a_1\dots a_r}\o j^rf$ have holomorphic extensions
to a neighborhood of $0$ and
$\Bigl(f^n_{a_1\dots a_r,\dots, a_1\dots a_r}\o j^rf\Bigr)(0)\ne 0$.
\endproclaim

\demo{Proof} We may suppose that $f(0)=0$ and $f\ne 0$. Let
$F:\Bbb C^p\to V$ be a local lift of $f$ at $0$. Since $F\ne 0$
there is a multi-index $A=(a_1,\dots,a_r)$
such that $\p_AF(0)\ne 0$. Then the functions
$f^1_{a_1\dots a_r}\o j^rf=F^{n-1}\p_AF$ and
$f^n_{a_1\dots a_r,\dots, a_1\dots a_r}\o j^rf=(\p_AF)^n$ are holomorphic
near $0$ and $f^n_{a_1\dots a_r,\dots, a_1\dots a_r}\o j^rf\ne 0$.

Suppose the conditions of the theorem are satisfied. By \thetag{\nmb|{1}},
the equality
$$
f^{n-1}f^n_{a_1\dots a_r,\dots, a_1\dots a_r}\o j^rf=
(f^1_{a_1\dots a_r}\o j^rf)^n \tag{\nmb:{2}}
$$
is satisfied in a neighborhood of $0$ outside the complex analytic set $f=0$.
Thus it is satisfied near $0$.
By assumption, the germ $f^n_{a_1\dots a_r,\dots, a_1\dots a_r}\o j^rf$
is invertible and the germs $f$ and $f^1_{a_1\dots a_r}\o j^rf$ are not invertible
in the ring $\frak F_{\Bbb C^p,0}$. Since this ring is factorial (see, for
example, \cit!{3}), \thetag{\nmb|{2}} implies that
$f^1_{a_1\dots a_r}\o j^rf$ divides $f$ in this ring. Then, by \thetag{\nmb|{2}},
$$
F:=\frac{f\root n\of{f^n_{a_1\dots a_r,\dots,a_1\dots a_r}\o j^rf}}
{f^1_{a_1\dots a_r}\o j^rf} \tag{\nmb:{3}}
$$
is a germ of a holomorphic function at $0$ and a local lift of $f$ at $0$.
By definition, $F$ is defined up to multiplication by
$\on{exp}\frac{2\pi i}{n}$.

Note that for $n=2$ we have, instead of \thetag{\nmb|{3}}, the following
simpler formula for the local lift $F$:
$$
F:=\frac{f^1_{a_1\dots a_r}\o j^rf}
{\root n\of{f^n_{a_1\dots a_r,\dots,a_1\dots a_r}\o j^rf}}.\qed
$$
\enddemo
We leave it to the reader to formulate the corresponding results for
formal and regular lifts.

\subhead{\nmb.{6.4}. Second Example}\endsubhead
Now we consider the dihedral groups.
Since the computations in this case are more complicated, we treat
only the tensor fields which are used in the conditions of the first order.

Let $G=\frak D_l$ $(l\geq 3)$ be the dihedral group acting on the real
Euclidean plane $\Bbb R^2$. We consider the corresponding complexification
of this action on $V=\Bbb C^2$.

Let $v=(x,y)\in\Bbb C^2$. We put $z:=x+iy$, $\widehat z:=x-iy$, and for
each integer $k>0$
$$
\on{R}z^k=\sum_{j=0}^{\left[\frac{k}{2}\right]}(-1)^j\binom{k}{2j}x^{k-2j}
y^{2j}\quad\text{and}\quad\on{I}z^k=\sum_{j=0}^{\left[\frac{k}{2}\right]}
(-1)^j\binom{k}{2j+1}x^{k-2j-1}y^{2j+1}.
$$
In particular, we have $\on{R}z=x$ and $\on{I}z=y$.
It is easily checked that $z^k=\on{R}z^k+i\on{I}z^k$ and
$\widehat z^k=\on{R}z^k-i\on{I}z^k$.

Choose the generators $\si_1$ and $\si_2$
of the ring $\Bbb C[V]^{\frak D_l}$ as follows:
$\si_1=\frac12z\widehat z=\frac12(x^2+y^2)$ and $\si_2=\frac1l\on{R}z^l=
\frac1{2l}(z^l+\widehat z^l)$.

By definition, $y_1=\si_1$ and $y_2=\si_2$ are the invariant
coordinates on $V$ which we will consider also as the coordinates on the
orbit space $Z=V/G$.

Denote by $\si_{1,s}$ and $\si_{2,s}$ the tensor fields $\ta_s$ from
\nmb!{4.3} for $\ta=\si_1$ and $\ta=\si_2$, respectively. To write
the conditions of the first order it suffices to calculate the tensor
fields $\si_{1,2}$ and $\si_{2,s}$ for $s=2,\dots,l$.

It is easily checked that for the coordinates $z$ and $\widehat z$ on $V$
we have
$$\aligned
J=-\frac{i}2\on{I}z^l\quad
S_1=\frac14(z_1\widehat z_2+z_2\widehat z_1),\quad
S_2=\frac1{2l}(z_1\cdots z_l+\widehat z_1\cdots\widehat z_l),\\
dv=-\frac1{i\on{I}z^l}\Bigl((\widehat z^{l-1}dy_1-zdy_2)\frac{\p}{\p z}+
(-z^{l-1}dy_1+\widehat zdy_2)\frac{\p}{\p\widehat z}\Bigr).
\endaligned\tag{\nmb:{1}}$$
One can put $\De:=(\on{I}z^l)^2=2^ly_1^l-l^2y_2^2$.

Furthermore we denote by $(dy_1)^p(dy_2)^q$ the symmetrized tensor product
of $p$ factors which are equal to $dy_1$ and $q$ factors which are
equal to $dy_2$.

Using \thetag{\nmb|{1}} we get
$$
S_1(dv,dv)= \frac1{(\on{I} z^l)^2}\Bigl(2^{l-1}y_1^{l-1}dy_1^2-2ly_2dy_1dy_2
+2y_1dy_2^2\Bigr).
$$
By the definition of $\si_{1,2}$, we have
$$
\si_{1,2}=\frac1{\tilde\De}\Bigl(2^{l-1}y_1^{l-1}dy_1^2-2ly_2dy_1dy_2
+2y_1dy_2^2\Bigr).
$$
Similarly, using \thetag{\nmb|{1}} we get
$$\align
&S_2(\underbrace{dv,\dots,dv}_{\text {$s$ times}},
\underbrace{v,\dots,v}_{\text {$l-s$ times}})=
\\&=
\frac{i^s}{2l(\on{I}z^l)^s}\sum_{t=0}^s\binom{s}{t}
\Bigl((-1)^{s-t}\widehat z^{(l-1)t}z^{l-t}+(-1)^tz^{(l-1)t}\widehat z^{l-t}\Bigr)
dy_1^tdy_2^{s-t}=
\\&=
\frac{(-1)^si^s}{2l(\on{I}z^l)^s}\Biggl((z^l+(-1)^s\widehat z^l)dy_2^s+
\\&\qquad\qquad\qquad\qquad+
\sum_{t=1}^s(-1)^{t}2^{l-t}\binom{s}{t}y_1^{l-t}
\Bigl(\widehat z^{l(t-1)}+(-1)^sz^{l(t-1)}\Bigr)dy_1^tdy_2^{s-t}\Biggr).
\endalign$$
Using the equality
$$
\widehat z^{l(t-1)}+(-1)^sz^{l(t-1)}=(\on{R}z^l-i\on{I}z^l)^{t-1}+
(-1)^s(\on{R}z^l+i\on{I}z^l)^{t-1}
$$
we get
$$
\multline
S_2(\underbrace{dv,\dots,dv}_{\text {$s$ times}},
\underbrace{v,\dots,v}_{\text {$l-s$ times}})
=
\frac{(-1)^si^s}{2l(\on{I}z^l)^s}\Biggl((z^l+(-1)^s\widehat z^l)dy_2^s+
\\
+\sum_{t=1}^s(-1)^{t}2^{l-t}\binom{s}{t}y_1^{l-t}\sum_{j=0}^{t-1}
\frac1{i^j}\binom{t-1}{j}(\on{R}z^l)^{t-j-1}(\on{I}z^l)^j(1+(-1)^{s+j})
dy_1^tdy_2^{s-t}\Biggr).
\endmultline\tag{\nmb:{2}}$$

Let $s=2r$. Then \thetag{\nmb|{2}} implies
$$\multline
\si_{2,2r}=\frac{(-1)^r}{\tilde\De^r}\Biggl(y_2dy_2^{2r}+\\
\sum_{t=1}^{2r}(-1)^t2^{l-t}
\binom{2r}{t}y_1^{l-t}\sum_{u=0}^{\left[\frac{t-1}{2}\right]}(-1)^u
l^{t-2u-2}\binom{t-1}{2u}\tilde\De^uy_2^{t-2u-1}dy_1^tdy_2^{2r-t}\Biggr).
\endmultline$$

Let $s=2r+1$.
Then \thetag{\nmb|{2}} implies
$$\multline
\si_{2,2r+1}=\frac{(-1)^r}{l\tilde\De^r}\Biggl(dy_2^{2r+1}+\\
\sum_{t=2}^{2r+1}(-1)^{t}2^{l-t}
\binom{2r+1}{t}y_1^{l-t}\sum_{u=1}^{\left[\frac{t}{2}\right]}
(-1)^ul^{t-2u}\binom{t-1}{2u-1}\tilde\De^{u-1}y_2^{t-2u}
dy_1^tdy_2^{2r+1-t}\Biggr).
\endmultline$$

By \nmb!{4.3} we can use the above tensor fields $\si_{1,2}$ and
$\si_{2,s}$ to write down the conditions of the first order for lifting for
the dihedral groups $\frak D_l$.

\remark{Remark} The conditions of lifting for the $\frak D_l$-module
$\Bbb C^2$ could be reduced to those of the first example \nmb!{6.3}. This follows from
the following formulas:
$$
(\on{I}z^l)^2=\tilde\De,\quad z^l=\on{R}z^l+i\on{I}z^l,\quad
\widehat z^l=\on{R}z^l-i\on{I}z^l.
$$
\endremark

\enddocument